\documentclass[12pt,a4paper]{article}

\usepackage[utf8]{inputenc}
\usepackage[T1]{fontenc}
\usepackage{epsfig}

\usepackage[margin=3cm]{geometry}

\usepackage{euler,times}
\usepackage{mathtools}

\newcommand{\footnoteremember}[2]{
\footnote{#2}
  \newcounter{#1}
  \setcounter{#1}{\value{footnote}}
}
\newcommand{\footnoterecall}[1]{
\footnotemark[\value{#1}]
}

\usepackage{float}
\usepackage{multirow}
\usepackage{pdflscape} 
\usepackage{booktabs}
\usepackage{xcolor,colortbl}
\usepackage{siunitx}
\usepackage{caption}
\usepackage{subcaption}
\usepackage{hyperref}

\usepackage{ amssymb }
\usepackage{amsmath}
\usepackage{pifont}% http://ctan.org/pkg/pifont

\usepackage[comma,authoryear]{natbib}
\usepackage{tikz}
\usetikzlibrary{plotmarks}
\usetikzlibrary{arrows}

\newcommand*\justify{%
  \fontdimen2\font=0.4em% interword space
  \fontdimen3\font=0.2em% interword stretch
  \fontdimen4\font=0.1em% interword shrink
  \fontdimen7\font=0.1em% extra space
  \hyphenchar\font=`\-% allowing hyphenation
}

\usepackage{algorithm}
\usepackage[noend]{algpseudocode}

\usepackage[toc,acronym,nonumberlist,nogroupskip]{glossaries}

\makeglossaries
\newacronym{bfgs}{BFGS}{Broyden-\-Fletcher-\-Goldfarb-\-Shanno}
\newacronym{cnn}{CNN}{Convolutional Neural Network}
\newacronym{dcm}{DCM}{Discrete Choice Model}
\newacronym{ffnn}{FFNN}{Feed-Forward Neural Network}
\newacronym{ml}{ML}{Machine Learning}
\newacronym{sgd}{SGD}{Stochastic Gradient Descent}
\newacronym[plural=VoT,firstplural=Values of Time (VoT)]{vot}{VoT}{Value of Time}
\newacronym{wma}{WMA}{Window Moving Average}
\newacronym{amabs}{AMABS}{Adaptive Moving Average Batch Size}
\newacronym{hamabs}{HAMABS}{Hybrid Adaptive Moving Average Batch Size}
\newacronym{lpmc}{LPMC}{London Passenger Mode Choice}
\newacronym{mtmc}{MTMC}{Mobility and Transport MicroCensus}

\title{
\vspace{-3cm}\epsfig{figure=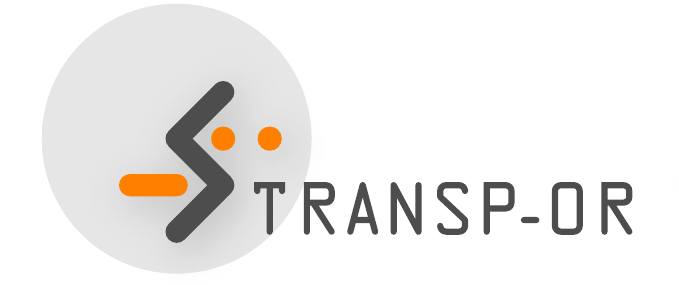,height=1.5cm} \hfill \epsfig{figure=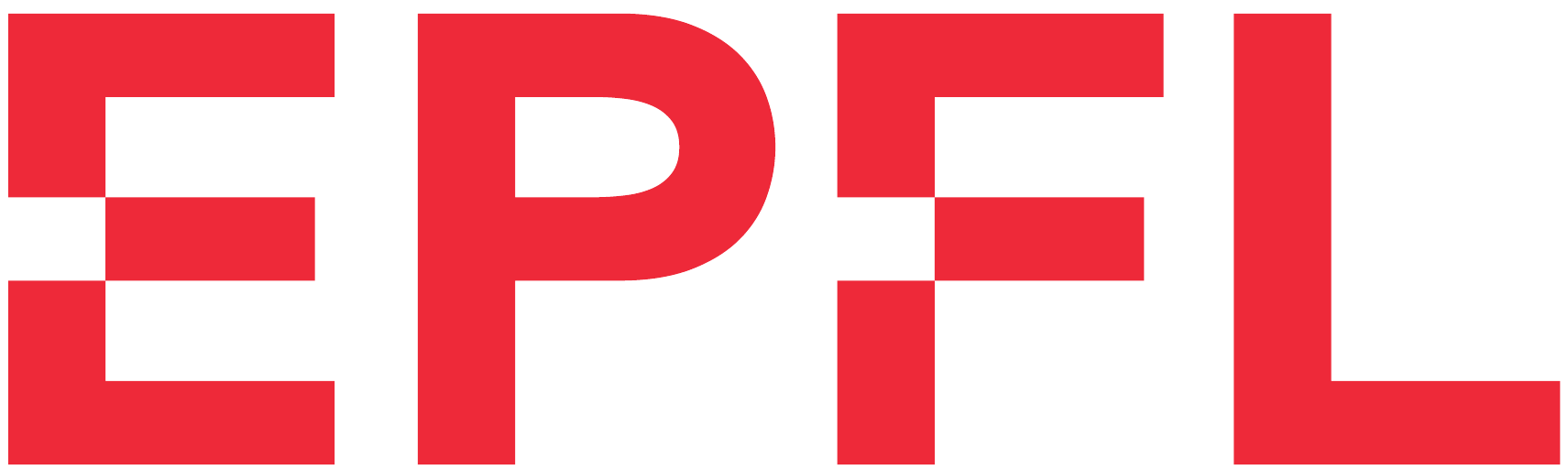,height=1.3cm}   \\*[-0.5cm] \mbox{}\hrulefill\mbox{} \\*[3cm] Estimation of Discrete Choice Models with Hybrid Stochastic Adaptive Batch Size Algorithms}
\author{
Gael Lederrey\footnoteremember{transp-or}{\'Ecole Polytechnique F\'ed\'erale de Lausanne (EPFL), School of Architecture, Civil and Environmental Engineering (ENAC), Transport and Mobility Laboratory,  Switzerland, \{gael.lederrey,tim.hillel,michel.bierlaire\}@epfl.ch}
\and
Virginie Lurkin\footnoteremember{TUE}{Eindhoven University of Technology, Department of Industrial Engineering \& Innovation Science, The Netherlands, v.j.c.lurkin@tue.nl}
\and
Tim Hillel\footnoterecall{transp-or}
\and
Michel Bierlaire\footnoterecall{transp-or}
}
\date{\today}

\begin{document}

\begin{titlepage}
\pagestyle{empty}

\maketitle
\vspace{2cm}

%%%%%%%%%%%%%%%%%%
%%%%%%%%%%%%%%%%%
%The report number is coded as YYMMDD where YY is the year, MM the
%month and DD the day.
%When a manuscript is finished, produce two version. One with the
%report number, and publish it as a technical report on our website,
%one without the report number, and submit it to a journal.
%In this case, just comment out the next lines.
\begin{center}
\small Report TRANSP-OR 191213 \\ Transport and Mobility Laboratory \\ School of Architecture, Civil and Environmental Engineering \\ Ecole Polytechnique F\'ed\'erale de Lausanne \\ \verb+http://transp-or.epfl.ch+
\end{center}

\thispagestyle{empty}

\end{titlepage}

%\tableofcontents

\newpage

\section*{Abstract}

The emergence of Big Data has enabled new research perspectives in the discrete choice community. While the techniques to estimate Machine Learning models on a massive amount of data are well established, these have not yet been fully explored for the estimation of statistical Discrete Choice Models based on the random utility framework. In this article, we provide new ways of dealing with large datasets in the context of Discrete Choice Models. We achieve this by proposing new efficient stochastic optimization algorithms and extensively testing them alongside existing approaches. We develop these algorithms based on three main contributions: the use of a stochastic Hessian, the modification of the batch size, and a change of optimization algorithm depending on the batch size. A comprehensive experimental comparison of fifteen optimization algorithms is conducted across ten benchmark Discrete Choice Model cases. The results indicate that the HAMABS algorithm, a hybrid adaptive batch size stochastic method, is the best performing algorithm across the optimization benchmarks. This algorithm speeds up the optimization time by a factor of 23 on the largest model compared to existing algorithms used in practice. The integration of the new algorithms in Discrete Choice Models estimation software will significantly reduce the time required for model estimation and therefore enable researchers and practitioners to explore new approaches for the specification of choice models. 

\paragraph{Keywords} Discrete Choice Models, Optimization, Stochasticity, Adaptive Batch Size, Hybridization
\section{Introduction}
\label{sec:intro}

The availability of more and more data for choice analysis is both a blessing and a curse. On the one hand, this data provides analysts with a great wealth of behavioral information. On the other hand, processing the increasingly large and complex datasets to estimate choice models presents new computational challenges. The \gls{ml} community has been thriving in dealing with vast amounts of data. It, therefore, seems natural to investigate \gls{ml} optimization algorithms to estimate choice models. 

The central ingredient of these algorithms is the \emph{stochastic gradient}. It consists in approximating the gradient of the log likelihood function (or any goodness of fit measure) using only a small subset of the data set. The gradient is said to be stochastic because the subset of data used to calculate it is drawn randomly from the full data set. A version of the steepest descent algorithm using this stochastic gradient is then applied to maximize the log likelihood function. Many variants have been proposed around this primary principle.

To illustrate the stochastic gradient, we consider applying stochastic gradient descent to a choice model with $J$ alternatives, and a data set of $N$ observations, each of them containing the vector of exploratory variables $x_n$ and the observed choice $i_n$. The choice model 
\begin{equation}
    P_n(i|x_n;\theta)
\end{equation}
provides the probability that individual $n$ chooses alternative $i$ in the context specified by $x_n$, where $\theta \in \mathbb{R}^K$ is a vector of $K$ unknown parameters, to be estimated from data. Typically, this is done using maximum likelihood estimation, where the log likelihood function $\mathcal{L}(\theta)$ is maximized:
\begin{equation}
\max_{\theta\in \mathbb{R}^K}  \mathcal{L}(\theta) = \max_{\theta\in \mathbb{R}^K} \sum_{n=1}^N \ln P_n(i|x_n;\theta).
\end{equation}
In the optimization literature, it is custom to define the algorithms for minimization problems. We follow the same convention, and consider the equivalent minimization problem: 
\begin{equation}
\label{eq:loglik}
\min_{\theta\in \mathbb{R}^K}  -\mathcal{L}(\theta)
\end{equation}
The gradient of the log likelihood function is 
\begin{equation}
    \nabla \mathcal{L}(\theta) = \sum_{n=1}^N \nabla \ln P_n(i|x_n;\theta).
\end{equation}
To obtain its stochastic version, we draw randomly, without replacement, a subset of $N'$ observations from the data that we call a \emph{batch}, and calculate:
\begin{equation}
    \label{eq:stoch_grad}
    \nabla_{N'} \mathcal{L}(\theta) = \sum_{n\in N'} \nabla \ln P_n(i|x_n;\theta).
\end{equation}

As shown in the above analysis, the variants of the stochastic gradient methods that are successful in \gls{ml} could be used as such to estimate the parameters of choice models. This could be used to decrease the time and computational cost of estimating choice models on large datasets. However, there are three critical differences between the two contexts which must be considered. Firstly, the parameters in choice models are a substantial output, and are used to estimate behavioral indicators such as \glspl{vot} and elasticities for the population. Conversely, parameters in \gls{ml} models typically have no behavioral interpretation, and only the model predictions are treated as a modeling output. It is, therefore, typical to allow choice model parameter estimates to converge during estimation to obtain the highest accuracy and precision of each individual parameter estimate. Conversely, in \gls{ml}, it is typical to restrict the model from converging fully to prevent overfitting, for example, by restricting the number of training epochs (one complete presentation of the data set). 

Secondly, choice models tend to have far fewer parameters than used in \gls{ml}. Complex choice models have hundreds of parameters, while the neural networks used in deep learning can involve millions of unknown parameters. 
There have been recent efforts to reduce the number of parameters in \gls{ml} models. For example, \cite{wu_prodsumnet:_2019} investigates simplifying \glspl{cnn}. The author can reduce the number of parameters using matrix decompositions from three million to just above three thousand, with only a small loss of precision. Nonetheless, these models are still complex and exceed the usual number of parameters used in choice models. 

Finally, choice data is typically collected using specific sampling strategies and designs of experiments. The objective is to obtain a representative sample while avoiding redundancies. Furthermore, the analyst usually has a model in mind when designing the data collection. In contrast, \gls{ml} techniques are often applied to datasets collected automatically or that were originally collected for other purposes (e.g., trip records from contactless payment cards), to detect patterns which have not previously been considered. It is, therefore, common to have a great deal of redundancy in the data. 
    
In this paper, we introduce a new algorithmic framework for the estimation of choice models, which addresses the key differences between the choice modeling and \gls{ml} contexts. Our framework includes three primary contributions:
\begin{itemize}
    \item[-] the use of a stochastic Hessian (that is, second derivative matrix),  which is possible thanks to the relatively low number of unknown parameters in discrete choice models,
    \item[-] the possible modification of the batch size from iteration to iteration, which is used to allow the high accuracy and precision in individual parameter estimates required in choice models,
    \item[-] a change of optimization algorithm depending on the size of the batch, which aims at finding the best trade-off between accuracy and efficiency.   
\end{itemize}

The rest of the paper is laid out as follows. In the next section, we describe optimization algorithms used both for the estimation of \glspl{dcm} and in \gls{ml}. Then, in Section~\ref{sec:methodology}, we present in detail the three ideas mentioned above and propose a catalog of variants of optimization algorithms for the estimation of choice models. In Section~\ref{sec:results}, we evaluate the performance of a series of algorithms on various choice models with large data sets. Finally, in Section~\ref{sec:conclusion}, we conclude this article and mention some further ideas that will be investigated in the future.
\renewcommand*{\arraystretch}{1} 
{\printglossary[type=\acronymtype, title=List of acronyms]}
\renewcommand*{\arraystretch}{1} 
\section{Literature Review}
\label{sec:lit_rev}

To the best of our knowledge, the maximum likelihood estimation of the parameters of choice models exclusively relies on deterministic algorithms that are variants of the line search and trust-region methods. In particular, the main estimation packages written in Python (Pandas Biogeme~\citep{bierlaire_biogeme:_2003, bierlaire_pandasbiogeme:_2018}, PyLogit~\citep{brathwaite_machine_2017}, and Larch~\citep{newman_computational_2018}) all use the \texttt{minimize} function from the package Scipy~\citep{jones_scipy:_2014}. This makes use of the quasi-Newton \gls{bfgs} algorithm. While this performs well for estimating small-to-medium-sized models, it struggles with larger, more complex models with more parameters. With the availability of larger and larger data sets, the performance of these standard methods is completely dominated by the time it takes to calculate the log likelihood function, its gradient, and its possible second derivative matrix. For example, \citet{hillel_understanding_2019} shows that fitting a \gls{ffnn} in the Tensorflow Python library is up to 200 times faster than estimating a Nested Logit model in Pandas Biogeme on the same data set containing 81'086 observations, despite the former having far more parameters. 

To understand how we may be able to estimate choice models more quickly on large datasets, we can look for inspiration from \gls{ml}. Datasets used in \gls{ml} (and in particular in Computer Vision) can contain millions of observations. For example, ImageNet, a collection of images labeled by hand, contains around fifteen million images. The sheer size of the full dataset ($\approx$150GB) prevents it from being stored in memory. As such, analyzing the full dataset is a significant computational challenge. As the full dataset cannot be stored in memory, the data must be batch processed. It explains the importance of the stochastic approach (calculate the gradient on a batch of data at each iteration) for \gls{ml}.To better understand the existing approaches used to estimate both choice models and \gls{ml}, we conduct a literature review of existing stochastic algorithms in Section~\ref{sec:existing_studies}. Then, in Section~\ref{sec:gaps}, we identify and discuss the gaps found in the literature for the optimization of choice models. We also link them to the three primary contributions presented in Section~\ref{sec:intro}.

\subsection{Overview of existing studies}
\label{sec:existing_studies}

In this section, we give a descriptive overview of the literature. We identify and review twenty-six studies that propose stochastic optimization algorithms. We focus specifically on quasi-Newton and second-order approaches while including examples of first-order algorithms. 
The selected studies, while not exhaustive, are believed by the authors to represent a broad overview of the existing optimization algorithms. 
For each of these algorithms, we analyze five aspects:
\begin{enumerate}
\item[-] the mathematical order, \emph{i.e.} if it uses a gradient-based (first-order), a quasi-Newton-based (1.5-th order), or a Newton-based step (second-order),
\item[-] if it uses an adaptive batch size technique or if the batch size is constant,
\item[-] if an experimental (numerical) assessment of the algorithm is conducted,
\item[-] a summary of the numerical applications of the algorithm.
\end{enumerate} 

The results of the review are displayed in tables \ref{tab:lit_rev_first} and \ref{tab:lit_rev_other}. Table \ref{tab:lit_rev_first} shows the details of the first-order stochastic algorithms, whilst Table~\ref{tab:lit_rev_other} shows the quasi-Newton (1.5th order) and second-order stochastic algorithms. Each table summarises the five considered features for each algorithm, alongside the reference and algorithm name. The following discussion evaluates the results in these tables. 

Of the twenty-six algorithms considered in this review, fourteen are first-order, six are quasi-Newton (1.5th order), and six are second order. We refer the reader to the article of \cite{ruder_overview_2016} for a precise overview of the first-order methods. While the review was focused on quasi-Newton and second-order algorithms, we find that the literature focuses predominantly on first-order approaches. We believe this focus is due to the speed of first-order algorithms for optimizing large, complex neural networks. 

Among the six second-order algorithms, three of them use the Hessian-Free (truncated Newton) optimization technique. This technique consists of approximating the problem using the second Taylor expansion and then solving it using a conjugate-gradient. The computation of the Hessian is approximated with a directional derivative and finite differences. The remaining three second-order stochastic methods use a subsampling method of the Hessian to avoid its heavy computation at each step. 

The majority of the algorithms use a fixed batch size, with only three out of the twenty-six algorithms using an adaptive batch size technique. All three articles which make use of an adaptive batch size are recent (2016-2018). None of the second-order methods and only one of the quasi-Newton algorithms make use of adaptive batch size. 

In terms of the analysis within the study, only five references do not provide numerical assessments of their algorithms. The first four references which do not include a numerical application are the oldest considered (1951-1992). The fifth reference, the algorithm RMSProp~\citep{tieleman_lecture_2012}, does not provide any theoretical nor numerical results since it was only presented in a lecture. This thus shows the importance of presenting numerical results. 

Finally, we can see that the algorithms in the study have been applied in multiple domains, but none explicitly for choice modeling. The earliest algorithms (4/26) were first applied to classical optimization problems. The remainder of the algorithms is applied to \gls{ml} problems, with the majority (13/26) being applied to neural networks. This shows the predominant focus on optimizing neural networks in the literature. 

\subsection{Gaps in knowledge}
\label{sec:gaps} 

As shown by the results of the literature review, the predominant focus of existing optimization research has been the optimization of neural networks, with none of the algorithms explicitly designed for the optimization of choice models. As discussed, there are substantial differences between the optimization of neural networks and \glspl{dcm}. In this section, we, therefore, assess the limitations of the existing algorithms in terms of optimizing \glspl{dcm}. Furthermore, we identify three gaps in knowledge in the existing research, which may enable higher-performing optimization algorithms for \glspl{dcm}. 

\subsubsection{Stochastic second order approaches}
\label{sec:stoch_hess_gaps}

\gls{ml} researchers have predominantly focused on first-order stochastic algorithms due to their speed when estimating parameters in large models. 
Existing research by \cite{lederrey_optimization_2018} shows that first-order stochastic methods are not able to achieve convergence on a logit model with ten parameters. The authors compare a gradient descent algorithm, a mini-batch \gls{sgd}, Adagrad~\citep{duchi_adaptive_2011}, and SAGA~\citep{defazio_saga:_2014}. They note that using a normalization on the parameters leads to more accurate results. However, they do not reach a sufficiently high precision in a reasonable amount of time. This failure is mainly due to the required precision for the convergence. As stated earlier, the parameter values are a key output of \glspl{dcm}. These parameters are used to calculate behavioral indicators such as the value-of-time (VOT) and elasticities. Therefore, it is critical to achieve convergence with high precision. 

Of the six second-order stochastic algorithms reviewed, none calculate the exact Hessian, with all using an approximation instead. 
However, the smaller number of parameters used in choice models compared to \gls{ml} models means computing the full Hessian for a batch of data is less computationally complex. There is, therefore, a need to investigate second-order stochastic approaches which compute the exact Hessian. The computation of the full Hessian could be used for choice models to obtain parameter estimates with the necessary precision for convergence. 

\subsubsection{Adaptive batch size}
\label{sec:abs_gaps}

None of the second-order methods found in the literature use an adaptive batch size. Furthermore, among the three algorithms proposing adaptive batch size methods, all of them couple the batch size with the \emph{learning rate} (a parameter specified in \gls{ml} estimation). While \cite{goyal_accurate_2018} shows that these two parameters are related, this coupling is specifically targeting \gls{ml} models, where the learning rate is often set to a fixed value or is slightly decreasing at each iteration.

The learning rate in machine learning is similar to the \emph{step size} used in classical optimization problems. Typically, in classical optimization, the step size is computed using more advanced techniques such as line search methods or trust-region methods. Since a high degree of precision is required for choice models, especially at the later stage of the optimization process, the step size should be separated from the batch size. Indeed, close to the optimum value, the optimization algorithm should both use a small step and as many data points as possible to achieve the highest precision. 

While line search and trust-region methods have already been applied to \gls{ml}, as demonstrated by \cite{rafati_trust-region_2018}, they have not yet been used in combination with adaptive batch size. Also, \cite{lederrey_optimization_2018, lederrey_snm:_2018} have demonstrated that the use of a full batch data is required at the end of the optimization process to achieve the appropriate precision for \glspl{dcm}. Therefore, it is required to develop an adaptive batch size technique with a second-order approach that does not interfere with the step size, and that eventually reaches the full dataset, as needed to achieve the required precision for choice models.

\subsubsection{Combined first and second order approaches}
\label{sec:hybrid_gaps}

All the algorithms presented in the review use the same algorithm throughout the optimization process. This can hinder the performance of the optimization since the requirements change during the process. Indeed, an optimization can be faster and less precise at the beginning of the optimization and then become more precise (and slower) close to the optimal solution. This can be partially achieved by using an adaptive batch size technique. However, this could be further addressed by switching the optimization algorithm at the right time to cope with the increased complexity due to a larger batch size. There is, therefore, a need to investigate hybrid algorithms, which switch between first and second-order approaches at the appropriate point in the optimization process. 

\begin{landscape}
    \vspace*{\fill}
\begin{table}[H]
    \centering
    \renewcommand\arraystretch{1}
    \begin{tabular}{l||c|c|c|c|m{6cm}}
        Reference & Name & Order & ABS & Experimental & Application \\ \hline \hline 
        \cite{robbins_stochastic_1951} & \gls{sgd} & 1st & & & Regression functions \\ \arrayrulecolor{gray}\hline
        \cite{polyak_methods_1964} & Momentum & 1st & & & Classical optimization problems \\ \arrayrulecolor{gray}\hline
        \cite{nesterov_method_1983} & NAG & 1st & & & Convex optimization problems \\ \arrayrulecolor{gray}\hline
        \cite{polyak_acceleration_1992} & Averaging & 1st & & & Classical optimization problems \\ \arrayrulecolor{gray}\hline
        \cite{duchi_adaptive_2011} & Adagrad & 1st & & \checkmark  & Image, text, and handwritten digit classification \\ \arrayrulecolor{gray}\hline
        \cite{zeiler_adadelta:_2012} & Adadelta & 1st & & \checkmark & Handwritten digit classification \\\arrayrulecolor{gray}\hline
        \cite{tieleman_lecture_2012} & RMSProp & 1st & & & None, first shown in a lecture \\\arrayrulecolor{gray}\hline
        \cite{schmidt_minimizing_2013} & SAG & 1st & & \checkmark & Binary classification \\\arrayrulecolor{gray}\hline
        \cite{defazio_saga:_2014} & SAGA & 1st & & \checkmark & Handwritten digit, binary, and multivariate classification \\\arrayrulecolor{gray}\hline
        \cite{kingma_adam:_2014} & Adam \& AdaMax & 1st & & \checkmark & Logistic Regression, Neural Networks, and Convolutional Neural Networks \\\arrayrulecolor{gray}\hline
        \cite{dozat_incorporating_2016} & Nadam & 1st & & \checkmark & Handwritten digit classification \\\arrayrulecolor{gray}\hline
           \cite{balles_coupling_2016} & CABS & 1st & \checkmark & \checkmark & Convolutional Neural Networks \\\arrayrulecolor{gray}\hline
   \cite{devarakonda_adabatch:_2017} & AdaBatch & 1st & \checkmark & \checkmark & Multiple Neural Networks architecture \\\arrayrulecolor{gray}\hline
        \cite{reddi_convergence_2018} & AMSGrad \& AdamNc & 1st & & \checkmark & Logistic Regression and Neural Networks \\
    \end{tabular}
    \caption{\label{tab:lit_rev_first} Analysis of first-order stochastic optimization algorithms included in the review. The term ``ABS'' stands for Adaptive Batch Size and ``Experimental'' for the experimental assessments.}
\end{table}
\end{landscape}

\begin{landscape}
    \vspace*{\fill}
\begin{table}[H]
    \centering
    \renewcommand\arraystretch{1.1}
    \begin{tabular}{l||c|c|c|c|m{6cm}}
        Reference & Name & Order & ABS & Experimental & Application \\ \hline \hline 
        \cite{bordes_erratum:_2010} & SGDQN & 1.5th & & \checkmark & Dense and sparse large datasets (PASCAL Large Scale challenge) \\\arrayrulecolor{gray}\hline
        \cite{martens_deep_2010} & HF & 2nd & & \checkmark & Image classification via Deep Learning models \\\arrayrulecolor{gray}\hline
        \cite{kiros_training_2013} & Stochastic HF & 2nd & & \checkmark & Classification and deep autoencoder tasks \\\arrayrulecolor{gray}\hline
        \cite{wang_stochastic_2014} & SQN \& RSQN & 1.5th & & \checkmark & Convex and non-convex optimization problems \\\arrayrulecolor{gray}\hline
        \cite{mokhtari_res:_2014} & RES & 1.5th & & \checkmark & Well- and ill-conditioned problems with large scale datasets \\\arrayrulecolor{gray}\hline
        \cite{you_investigation_2014} & SHF & 2nd & & \checkmark & Speech recognition \\\arrayrulecolor{gray}\hline
        \cite{keskar_adaqn:_2016} & adaQN & 1.5th & & \checkmark & Recurrent Neural Networks \\\arrayrulecolor{gray}\hline
        \cite{agarwal_second-order_2016} & LiSSA & 2nd & & \checkmark & Handwritten digit and multivariate classification \\\arrayrulecolor{gray}\hline 
        \cite{mutny_stochastic_2016} & ISSA & 2nd & & \checkmark & Least square estimators \\\arrayrulecolor{gray}\hline 
        \cite{ye_nestrovs_2017} & AccRegSN & 2nd & & \checkmark & Least square regressions \\\arrayrulecolor{gray}\hline 
        \cite{gower_accelerated_2018} & hBFGS & 1.5th & & \checkmark & Matrix Inversion problems \\\arrayrulecolor{gray}\hline
        \cite{bollapragada_progressive_2018} & PBQN & 1.5th & \checkmark & \checkmark & Logistic Regressions and Neural Networks \\
    \end{tabular}
    \caption{\label{tab:lit_rev_other} Analysis of  1.5th- and second-order stochastic optimization algorithms included in the review. The term ``ABS'' stands for Adaptive Batch Size and ``Experimental'' for the quantitative assessments.}
\end{table}
    \vspace*{\fill}
\end{landscape}

\subsection{Summary}

In the previous section, we identify gaps in knowledge in the literature across three themes; second-order stochastic approaches, adaptive batch size, and combined first and second-order approaches. Researchers have already investigated solutions for some of these limitations individually. However, we could not find any research investigating their combination in a systematic approach. We thus aim at combining the three primary contributions stated in the introduction, Section~\ref{sec:intro}, in a new algorithm for the optimization of choice models. Thus, this algorithm should incorporate the following concepts:
\begin{itemize}
    \item[-] the use of a second-order stochastic approach which calculates the true Hessian for each batch,
    \item[-] the use of adaptive batch size for second-order approaches which do not couple the batch size to the learning rate,
    \item[-] the use of a hybrid approach that combines first and second-order algorithms and switches between them at the appropriate time.
\end{itemize}
\section{Methodology}
\label{sec:methodology}

This section introduces our novel algorithmic framework for estimating \glspl{dcm}. Sections~\ref{sec:ls} and \ref{sec:tr} provide a reminder of the two standard optimization techniques: \emph{line search} and \emph{trust regions}. Section~\ref{sec:methodo_approach} shows how we construct the new hybrid stochastic algorithms with adaptive batch size. 

\subsection{Line search methods}
\label{sec:ls}

Line search optimization methods combine a descent direction with a line search. The iterates are $\theta_{k+1} = \theta_k + \alpha_k d_k$, where $d_k$ is a descent direction obtained by preconditioning the gradient. $d_k$ is thus defined as:
\begin{equation}
\label{eq:precond}
d_k = - D_k \nabla \mathcal{L}(\theta_k)
\end{equation}
where $D_k$ is a positive definite matrix.

There are typically three ways to select $D_k$:
\begin{itemize}
\item[-]  {\bf Steepest descent} methods define $D_k$ as the identity matrix. In that case, the descent direction is the opposite of the gradient. 
\item[-] {\bf Newton's method} assumes that $D_k$ is the inverse of the second derivative matrix (possibly perturbed to make it definite positive). 
\item[-] Quasi-Newton methods assume that $D_k$ is a secant approximation of (the inverse of) the second derivative matrix, updated at each iteration. Among the many secant methods, we consider the \textbf{\gls{bfgs}} algorithm for which, according to~\citet{fletcher_practical_1987}, the approximation is given by:
\begin{equation}
    \label{eq:bfgs}
    B_{k+1} = B_k + \frac{y_k y_k^T}{y_k^T s_k} - \frac{B_k s_k s_k^T B_k^T}{s_k^T B_k, s_k},
\end{equation}
with $s_k = \theta_{k+1} - \theta_k$ and $y_k = \nabla \mathcal{L}(\theta_{k+1}) - \nabla \mathcal{L}(\theta_k)$. 

A slightly different version of \gls{bfgs} consists in approximating the inverse of the Hessian. In that case, the name {\bf \gls{bfgs}$^{-1}$} is used and the approximation is given by: 
\begin{equation}
    \label{eq:bfgs-inv}
    B_{k+1}^{-1} = B_k^{-1} + \frac{\left(s_k^Ty_k + y_k^TB_k^{-1}y_k\right)\left(s_ks_k^T\right)}{\left(s_k^Ty_k\right)^2} - \frac{B_k^{-1}y_ks_k^T + s_ky_k^TB_k^{-1}}{s_k^Ty_k}.
\end{equation}
\end{itemize}

The step is calculated with an inexact line search method, based on the two Wolfe conditions~\citep{wolfe_convergence_1969,wolfe_convergence_1971}. The first condition guarantees that the step gives a sufficient decrease in the objective function, while the second one makes sure that unacceptably small steps are ruled out. 

\subsection{Trust-region methods}
\label{sec:tr}
Trust-region methods define a region around the current search point, where a quadratic approximation of the function value is "trusted" to be correct, and steps are chosen to optimize the function model within this region. There exist multiple ways to compute the quadratic approximation of the function value. The most common choice is a quadratic function of the type:
\begin{equation}
    m_k(\theta_k + d_k) = \mathcal{L}(\theta_k) + \nabla \mathcal{L}(\theta_k)^Td_k + \frac{1}{2} d_k^T B_k d_k
\end{equation}
where $B_k$ is either the Hessian $\nabla^2 \mathcal{L}(\theta_k)$ or an approximation of it. If the Hessian is chosen, the name \textbf{trust-region} algorithm is used. If an approximation of the Hessian is used it is referred to as a \textbf{quasi-Newton trust-region}. 

The size of the trust region is modified during the search, based on how well the quadratic model agrees with the actual objective function value. Following \citet{conn_trust_2000}, the region is modified conditional to the ratio~$\rho_k$ of the actual function reduction to the reduction predicted by the quadratic model:
\begin{equation}
    \rho_k = \frac{\mathcal{L}(\theta_k) - \mathcal{L}(\theta_k + d_k)}{m_k(\theta_k) - m_k(\theta_k + d_k)}.
\end{equation}
Given the ratio~$\rho_k$, the decision to change the trust region is based on the following rules:
\begin{equation}
    \Delta_k = \begin{cases}
                \gamma_1 \Delta_k & \text{if } \rho_k \geq \eta_2,\\
                \Delta_k & \text{if } \rho_k \in [\eta_1, \eta_2),\\
                \gamma_2 \Delta_k & \text{if } \rho_k < \eta_1,\\
                \end{cases}
\end{equation}
where $\gamma_1, \gamma_2, \eta_1$ and $\eta_2$ are all \emph{a priori} defined parameters. 

Intuitively, these rules imply that when the quadratic model is a good predictor of the function value (ratio close to 1), the region is not modified (second case). In contrast, when the ratio is too small (third case), respectively too large (first case), the quadratic model is no longer be a good predictor, and so the region is reduced, respectively extended. 

\subsection{Hybrid Stochastic Algorithms with Adaptive Batch Size}
\label{sec:methodo_approach}
We present our new algorithmic framework for optimization of choice models by discussing its three main contributions: 
\begin{itemize}
    \item[-] the use of a stochastic hessian,
    \item[-] the possible modification of the batch size from iteration to iteration,
    \item[-] the change of optimization algorithm depending on the size of the batch.   
\end{itemize}

\subsubsection{Stochastic Hessian}
\label{sec:stoch_hess_methodo}
Inspired by the technique of stochastic gradient, our algorithmic framework includes a stochastic Hessian, defined as:
\begin{equation}
    \label{eq:stoch_hess}
    \nabla_{N'}^2 \mathcal{L}(\theta) = \sum_{n\in N'} \nabla^2 \ln P_n(i|x_n;\theta),
\end{equation}
where $N'$ is a subset of observations, drawn randomly, without replacement, from the full dataset. 

A stochastic Hessian can be included in \gls{sgd} algorithms to create a {\bf Stochastic Newton Method} (SNM), as in~\cite{lederrey_snm:_2018}, or a {\bf Stochastic Trust-Region} (STR) algorithm.

\subsubsection{Adaptive Batch Size}
\label{sec:abs_methodo}

Intuitively, this technique increases the batch size when the algorithm is not able anymore to improve the value of the objective function. This generally happens for two reasons. Either a local optimum in the current neighborhood has been reached, or the stochastic nature of the gradient and Hessian precludes the algorithm from making further progress. 

We propose to modify the size of the batch as the algorithm proceeds with a \gls{wma} technique. \gls{wma} is a simple smoother, consisting of averaging the values of a time series across a window of $W$ consecutive observations, thereby generating a series of averages. To better capture change in the data, more importance is usually given to the most recent iterations. The \gls{wma} at the $k$-th iteration is given by:
\begin{equation}
    \label{eq:wma}
\text{WMA}_{k,W} = \begin{dcases}
\frac{\sum_{i=0}^{\mathcal{W}-1} (\mathcal{W}-i)\mathcal{L}(\theta_{k-i})}{\sum_{i=1}^{\mathcal{W}} i} & \text{if } k \geq \mathcal{W}, \\
\frac{\sum_{i=0}^{k-1} (k-i)\mathcal{L}(\theta_{k-i})}{\sum_{i=1}^{k} i} & \text{otherwise.}
\end{dcases}
\end{equation}
Note that for all first $W$ iterations, the size of the window is reduced to the iteration number, \emph{i.e.} $W = k$. 

The decision to increase the batch size is then based on a successive lack of progress in the past iterations. The progress at the $k$-th iteration is defined as:
\begin{equation}
    \label{eq:impr}
    \mathcal{I}_k = \frac{\text{WMA}_{k-1,W}-\text{WMA}_{k,W}}{\text{WMA}_{k-1,W}}.
\end{equation}
We consider that there is a lack of progress when ${I}_k$ is less than a threshold ${I}_k < \Delta$. After $C$ iterations with a lack of progress, the batch size is increased by a factor of $\tau$. 

We call the above algorithm \gls{amabs} and its pseudocode is given in Algorithm~\ref{algo:amabs}. The \gls{amabs} technique can be used in any stochastic optimization algorithm. For example, using the principle of stochastic Hessian, an {\bf \gls{amabs} Newton's method} or an {\bf \gls{amabs} Trust-Region} method can be defined. 

\begin{algorithm*}[h!]
\caption{Adaptive Moving Average Batch Size (\gls{amabs})}\label{algo:amabs}
\begin{algorithmic}[1]
\Require{
\begin{itemize}
\itemsep0em 
\item[]
\item[-] Current iteration index: $k$,
\item[-] Function value at iteration $k$: $\mathcal{L}(\theta_k)$,
\item[-] Current batch size: $N'_k$,
\item[-] Size of the full dataset: $N_{\max}$, 
\item[-] Size of the window: $W$ (default: $10$), 
\item[-] Threshold for successfull iterations: $\Delta$ (default: $1\%$), 
\item[-] Maximum number of unsuccessfull iterations with the same batch size: $C$ (default: $2$), 
\item[-] Expansion factor for the batch size: $\tau$ (default: $2$), 
\end{itemize}
}
\Ensure{New batch size: $N'_{k+1}$}
\State Set the counter at 0 ($c=0$) at the initialization of the algorithm
\Function{AMABS}{}
\State Compute WMA$_{m,W}$, as in Equation~\ref{eq:wma}.
\If {$k > 0$} \Comment{At least two iterations required}
\State Compute $\mathcal{I}_m$ as in Equation~\ref{eq:impr} using WMA$_{k,W}$ and WMA$_{k-1,W}$
    \If {$\mathcal{I}_M < \Delta$} \Comment{Count the consecutive steps under threshold.}
        \State $c = c+1$
    \Else
        \State $c = 0$ 
    \EndIf
    \If {$c == C$} \Comment{Update the batch size}
        \State $c = 0$
        \State $N'_{k+1} = \min(\tau\cdot N',N_{\max})$
    \Else
        \State $N'_{k+1} = N'_k$
    \EndIf
\EndIf
\State \textbf{return} $N'_{k+1}$
\EndFunction
\end{algorithmic}
\end{algorithm*}

\subsubsection{Hybridization}
\label{sec:hybrid_methodo}

The use of the \gls{amabs} technique naturally brings the opportunity of using different optimization algorithms based on the batch size. Indeed, at the beginning of the optimization process, when small batch size is used, it makes sense to use the exact second derivative matrix, as it is relative cheap to calculate. When the batch size becomes large, the calculation time for the exact Hessian is prohibitive, and methods relying on quasi-Newton approximation are preferred. It makes, therefore, sense to switch to a less precise but faster algorithm. Our hybrid algorithm determines the algorithm to use based on the percentage of data used in a batch. The pseudocode is shown in Algorithm~\ref{algo:hybrid}. 

\begin{algorithm*}[h!]
\caption{Hybrid algorithm with \gls{amabs} method}\label{algo:hybrid}
\begin{algorithmic}[1]
\Require{
\begin{itemize}
\itemsep0em 
\item[]
\item[-] Log likelihood function: $\mathcal{L}(\theta)$,
\item[-] Gradient of the log likelihood: $\nabla \mathcal{L}(\theta)$,
\item[-] Hessian of the log likelihood: $\nabla^2 \mathcal{L}(\theta)$,
\item[-] Initial parameter value: $\theta_0$,
\item[-] Algorithm generating a candidate for the next iteration using only the gradient (first-order or quasi-Newton method): \texttt{generateCandidateFirstOrder}$(\theta,\nabla \mathcal{L}(\theta))$,
\item[-] Algorithm generating a candidate for the next iteration using the gradient and the Hessian (second-order method): \texttt{generateCandidateSecondOrder}$(\theta,\nabla \mathcal{L}(\theta),\nabla^2 \mathcal{L}(\theta))$,
\item[-] Parameters specific to the algorithm:
\begin{itemize}
\itemsep0em 
\item[\textasteriskcentered] Initial batch size: $N_0'$ (default: $1000$),
\item[\textasteriskcentered] Size of the full dataset: $N_{\max}$, 
\item[\textasteriskcentered] Size of the window: $W$ (default: $10$), 
\item[\textasteriskcentered] Threshold for successfull iterations: $\Delta$ (default: $1\%$), 
\item[\textasteriskcentered] Maximum number of unsuccessfull iterations with the same batch size: $C$ (default: $2$), 
\item[\textasteriskcentered] Expansion factor for the batch size: $\tau$ (default: $2$), 
\item[\textasteriskcentered] Threshold for hybridization: $\Delta_H$ (default: $30\%$), 
\item[\textasteriskcentered] Threshold for stopping criterion: $\varepsilon$ (default: $10^{-6}$)
\end{itemize}
\end{itemize}
}
\Ensure{Optimized parameters: $\theta^*$}
\Function{Iteration}{}
\If {$N'_k/N_{\max} > \Delta_H$}
\State $\theta_{k+1} = \texttt{\justify generateCandidateFirstOrder}(\theta_k,\mathcal{L}(\theta), \nabla \mathcal{L}(\theta))$
\Else 
\State $\theta_{k+1} = \texttt{\justify generateCandidateSecondOrder}(\theta_k,\mathcal{L}(\theta), \nabla \mathcal{L}(\theta),\nabla^2 \mathcal{L}(\theta))$
\EndIf
\State $N_{k+1}' = AMABS(k, \mathcal{L}(\theta_k), N'_k, N_{\max}, W, \Delta, C, \tau)$
\State Stop the optimization if $\nabla_{\text{rel}}\mathcal{L}(\theta_{k+1}) < \varepsilon$. 
\EndFunction
\end{algorithmic}
\end{algorithm*}

The usual initial parameter value $\theta_0$ is an array of zeros. The last two parameters introduced in Algorithm~\ref{algo:hybrid}, $\Delta_H$ and $\varepsilon$, are further discussed in Section~\ref{sec:parameters_algo}. The stopping criterion $\nabla_{rel} \mathcal{L}(\theta)$ is discussed at the end of this section. The functions \texttt{\justify generateCandidateFirstOrder} (line 4) and \texttt{\justify generateCandidateSecondOrder} (line 7) return the parameter values for the next iteration. The generic pseudocode for these methods is given in Algorithm~\ref{algo:gen_cand}. The function \texttt{\justify generateCandidateFirstOrder} uses the function \texttt{\justify generateCandidate} with a first-order or a quasi-Newton method and thus does not use the Hessian ($\nabla^2 \mathcal{L}(\theta)$). The function \texttt{\justify generateCandidateSecondOrder} uses the function \texttt{\justify generateCandidate} with the Hessian. 

\begin{algorithm*}[h!]
\caption{Pseudocode for the function \texttt{generateCandidate}}\label{algo:gen_cand}
\begin{algorithmic}[1]
\Require{
\begin{itemize}
\itemsep0em 
\item[]
\item[-] Parameter value at iteration $k$: $\theta_k$
\item[-] Log likelihood function: $\mathcal{L}(\theta)$,
\item[-] Gradient of the log likelihood: $\nabla \mathcal{L}(\theta)$,
\item[-] Possibly, Hessian of the log likelihood: $\nabla^2 \mathcal{L}(\theta)$
\end{itemize}
}
\Ensure{Next parameter value: $\theta_{k+1}$}
\Function{\texttt{generateCandidate}}{}
\State Compute direction $d_k$ using $\nabla \mathcal{L}(\theta_k)$ and/or $\nabla^2 \mathcal{L}(\theta_k)$ as in Equation~\ref{eq:precond}.
\State Compute step size $\alpha_k$ using either a line search method, Section~\ref{sec:ls}, or a trust-region method, Section~\ref{sec:tr}.
\State Compute $\theta_{k+1} = \theta_k + \alpha_k d_k$
\EndFunction
\end{algorithmic}
\end{algorithm*}

\paragraph{Stopping Criterion}
Both line search and trust-region methods require a stopping criterion to determine when to stop iterating. A relative gradient-based stopping criterion ensures that the algorithm stops when the norm of the relative gradient is below some threshold, $\varepsilon$:

\begin{equation}
    \label{eq:stop_crit}
    ||\nabla_{rel} \mathcal{L}(\theta) || \leq \varepsilon.
\end{equation}
A sufficiently small value is chosen for $\varepsilon$, typically in the range $[10^{-6}, 10^{-8}]$. A more detailed discussion of the stopping criterion can be found in~\citet{dennis_numerical_1996} (see Chapter 7.2 page 159).

\subsection{Summary of algorithms}
\label{sec:sum_algos}

As depicted in Table~\ref{tab:algorithms}, the new algorithmic framework presented in Sections~\ref{sec:ls}, \ref{sec:tr}, and \ref{sec:hybrid_methodo} allows us to compare the performance of 15 different algorithms. These algorithms can be split into three categories:
\begin{itemize}
\item[-] \emph{standard non-stochastic algorithms} (first 6 algorithms) -- deterministic algorithms which are commonly used in the \gls{dcm} community and are therefore considered as benchmarks. For example, the current version of Biogeme uses the Python package Scipy~\citep{jones_scipy:_2014} with the \gls{bfgs}$^{-1}$ implementation. We will start with a comparison of the performance of these standard algorithms before moving to the comparison with stochastic approaches. 
\item[-] \emph{stochastic algorithms} (next 6 algorithms) -- algorithms based on the \gls{amabs} method presented in Section~\ref{sec:abs_methodo}. These methods are used to show the improvement in terms of estimation time over their non-stochastic counterparts. The algorithms \texttt{NM-ABS} and \texttt{TR-ABS} also make use of the stochastic Hessian, as presented in Section~\ref{sec:stoch_hess_methodo}. The added value of using stochastic Hessian will therefore also be discussed. 
\item[-] \emph{hybrid stochastic methods} (last 3 algorithms) -- algorithms which combine two separate optimization algorithms as presented in Section~\ref{sec:hybrid_methodo}. Three types of hybridization are investigated and discussed: (i)~Newton's method and \gls{bfgs}, (ii)~Trust-Region method and \gls{bfgs}, and (iii)~Newton's method and \gls{bfgs}$^{-1}$. The first algorithm always corresponds to the function \texttt{\justify generateCandidateSecondOrder} and the second to the function \texttt{\justify generateCandidateFirstOrder} in Algorithm~\ref{algo:hybrid}.
\end{itemize}

\begin{table}[h!]
    \centering
    \small
    \begin{tabular}{l||c|c|p{8.5cm}}
        Name & Order & \gls{amabs} & Description \\ 
        \hline
        \hline
        \texttt{GD} & 1st & & Steepest descent algorithm. \\
        \texttt{BFGS} & 1.5th & & \gls{bfgs} algorithm using Eq.~\ref{eq:bfgs}. \\
        \texttt{BFGS$^{-1}$} & 1.5th & & \gls{bfgs}$^{-1}$ algorithm using Eq.~\ref{eq:bfgs-inv}. \\
        \texttt{TR-BFGS} & 1.5th & & quasi-Newton trust-region method with \gls{bfgs} (Eq.~~\ref{eq:bfgs}). \\
        \texttt{NM} & 2nd & & Newton's method. \\
        \texttt{TR} & 2nd & & Trust-region method. \\ 
        \hline
        \hline 
        \texttt{GD-ABS} & 1st & \checkmark & Stochastic steepest descent with \gls{amabs}. \\
        \texttt{BFGS-ABS} & 1.5th & \checkmark & BFGS algorithm (Eq.~\ref{eq:bfgs}) with \gls{amabs}. \\
        \texttt{BFGS$^{-1}$-ABS} & 1.5th & \checkmark & BFGS$^{-1}$ algorithm (Eq.~\ref{eq:bfgs-inv}) with \gls{amabs}. \\
        \texttt{TR-BFGS-ABS} & 1.5th & \checkmark & Trust-Region with \gls{bfgs} (Eq.~\ref{eq:bfgs}) and \gls{amabs}. \\
        \texttt{NM-ABS} & 2nd & \checkmark & Newton with \gls{amabs}. \\
        \texttt{TR-ABS} & 2nd & \checkmark & Trust-region with \gls{amabs}. \\ \hline\hline
        \texttt{H-NM-ABS} & Hybrid & \checkmark & Hybridization: Newton + \gls{bfgs} (Eq.~\ref{eq:bfgs}). \\
        \texttt{H-TR-ABS} & Hybrid & \checkmark & Hybridization: trust-region + \gls{bfgs} (Eq.~\ref{eq:bfgs}). \\
        \texttt{HAMABS} & Hybrid & \checkmark & Hybridization: Newton + \gls{bfgs}$^{-1}$ (Eq.~\ref{eq:bfgs-inv}).
    \end{tabular}
    \normalsize
    \caption{\label{tab:algorithms}Overview of all algorithms used for the optimization of \glspl{dcm}. A small description of the algorithms is provided as well as their order and if it is including the adaptive batch size method.}
\end{table}
\section{Results}
\label{sec:results}
This section presents the results of our experiments. Before showing the numerical results, we start by explaining our experimental design, \textit{i.e.}, the algorithms and datasets used in our experiments, as well as the implementation details. 

\subsection{Experimental Design}
\label{sec:exp_design}

We collect empirical evidence of the behavior of the 15 algorithms, in Table~\ref{tab:algorithms}, by observing their performance on 10 different choice models, presented in Table~\ref{tab:models}. As shown in the column \textit{Data}, two different sources of data are used. The first 9 choice models are estimated on data obtained from the \gls{lpmc} dataset. This dataset, collected by \cite{hillel_recreating_2018}, contains mode choice on an urban multi-modal transport network, from April 2012 to March 2015. Three sub-datasets have been created to study the impact of the size of the dataset on the performance of the different algorithms: 
\begin{itemize}
    \item[-] a small dataset (S), that contains observations from year 2012 (27'478 observations),
 \item[-] a medium dataset (M), that contains observations from years 2012 to 2013 (54'766 observations),
\item[-]  a large dataset (L), that contains all observations, from year 2012 to 2015 (81'766 observations). 
\end{itemize}
In order to analyze the impact of the number of parameters on the estimation time, we compare three logit models from \cite{hillel_understanding_2019}:
\begin{itemize}
    \item[-] the \texttt{\gls{lpmc}\_DC} model that contains 13 parameters to be estimated,
 \item[-] the \texttt{\gls{lpmc}\_RR} model that contains 54 parameters to be estimated,
\item[-]  the \texttt{\gls{lpmc}\_Full} model that contains 100 parameters to be estimated. 
\end{itemize}
   
Finally, a tenth choice model was estimated on the \gls{mtmc} dataset, a statistical survey of the travel behavior of the Swiss population (\cite{danalet_mobility_2018}). We use the most recent version of the survey, collected in 2015. This model provides an exciting opportunity to study the efficiency of all algorithms presented in Table~\ref{tab:algorithms} for estimating a rather large choice model (almost 250 parameters), on a medium-size dataset (56'915 observations).

\begin{table}[h]
  \centering
  \renewcommand\arraystretch{1.2}
  \setlength{\tabcolsep}{0.5em}
  \begin{tabular}{l|ccccc}
   Names & \#Parameters & Data & \#Observations \\ \hline \hline
   \texttt{\gls{lpmc}\_DC\_S} & 13 & \gls{lpmc} & 27'478 \\
   \texttt{\gls{lpmc}\_DC\_M} & 13 & \gls{lpmc} & 54'766 \\
   \texttt{\gls{lpmc}\_DC\_L} & 13 & \gls{lpmc} & 81'086 \\ \hline
   \texttt{\gls{lpmc}\_RR\_S} & 54 & \gls{lpmc} & 27'478 \\
   \texttt{\gls{lpmc}\_RR\_M} & 54 & \gls{lpmc} & 54'766 \\
   \texttt{\gls{lpmc}\_RR\_L} & 54 & \gls{lpmc} & 81'086 \\ \hline
   \texttt{\gls{lpmc}\_Full\_S} & 100 & \gls{lpmc} & 27'478 \\
   \texttt{\gls{lpmc}\_Full\_M} & 100 & \gls{lpmc} & 54'766 \\
   \texttt{\gls{lpmc}\_Full\_L} & 100 & \gls{lpmc} & 81'086 \\ \hline \hline
   \texttt{MTMC} & 247 & \gls{mtmc} & 56'915 \\
  \end{tabular}
  \caption{\label{tab:models} Summary of the models used for the performance analysis. The number of parameters is provided as well as the dataset and the number of observations used in the model. All the models are logit models.}
\end{table}    

\subsection{Implementation details} 
All models are estimated on a single node in a supercomputer (18 Cores Skylake Processor@2.30 GHz, 192GB) for each algorithm. We include a stopping criterion on the maximum number of epochs (1,000 epochs). This is done to avoid extremely long computation time for algorithms that would struggle in achieving convergence for certain models. Also, since some of these algorithms are stochastic, the speed of convergence may differ on the same optimization task. Thus, each stochastic algorithm is used to optimize each model 20 times. We impose an upper limit on the execution time for the 20 estimation process: 12 hours for the models \texttt{\gls{lpmc}\_DC}, 24 hours for the models \texttt{\gls{lpmc}\_RR}, 36 hours for the models \texttt{\gls{lpmc}\_Full}, and 48 hours for the model \texttt{\gls{mtmc}}. Finally, since we want to be able to compare the efficiency of our algorithms with state-of-the-art \gls{dcm} software, we also optimize all the models in Table~\ref{tab:models} with Biogeme and Scipy within the same rules. All the results are presented and discussed in Section~\ref{sec:performance}. The code implementing all the algorithms in Table~\ref{tab:algorithms} can be found on Github at \href{https://github.com/glederrey/HAMABS}{https://github.com/glederrey/HAMABS}.

\subsection{Performance analysis}
\label{sec:performance}

In this section, we analyze the performance of the fifteen algorithms reported in Table~\ref{tab:algorithms}. For ease of comparison, we decided to use a graphical approach named \textit{performance profiles} to benchmark these algorithms. As stated by \cite{beiranvand_best_2017}, performance profiles are a great tool to analyze algorithms in terms of efficiency, robustness, and probability of successfully performing a required task. The concept of performance profile is presented in Section~\ref{sec:concep_perf_prof}, while results are analyzed in Section~\ref{sec:res_perf_prof}.

\subsubsection{Performance profiles}
\label{sec:concep_perf_prof}

Performance profiles are introduced by \cite{dolan_benchmarking_2002} and are now a recurring tool used to compare the performance of optimization algorithms. They are used to compare the performance of a set of optimization algorithms, $\mathcal{A}$, on a set of optimization problems, $\mathcal{P}$. For each pair $(p,a)\in \mathcal{P} \times \mathcal{A}$, they define a performance measure $t_{p,a} > 0$ whose large value indicates poor performance. Classical measures of performance are the execution time or the number of epochs. 

In addition, a convergence test $\mathcal{C}_{p,a}$ states if algorithm $a$ was able to optimize problem $p$. For each optimization problem $p$ and optimization algorithm $a$, the performance ratio is defined as
\begin{equation}
    \label{eq:perf_ratio}
    r_{p,a} = \begin{cases}
    \frac{t_{p,a}}{\min_{a\in\mathcal{A}} t_{p,a}} & \text{if $\mathcal{C}_{p,a}$ passed,} \\
    \infty & \text{if $\mathcal{C}_{p,a}$ failed.}
    \end{cases}
\end{equation}
This leads to $r_{p,a}=1$ for the best algorithm and $r_{p,a}=\infty$ for all algorithms $a$ unable to solve problem $p$. The performance profile of an algorithm $a$ is finally defined as
\begin{equation}
    \label{eq:perf_prof}
    \rho_a(\pi) = \frac{\left| p\in \mathcal{P}: r_{p,a} \leq \pi \right|}{|\mathcal{P}|} 
\end{equation}
where $|\mathcal{P}|$ is the cardinality of the set $\mathcal{P}$. $\rho_a(\pi)$ represents the proportion of problems for which the performance ratio $r_{p,a}$ for algorithm $a\in\mathcal{A}$ is within a factor $\pi\in\mathbb{R}$ of the best possible performance ratio.

Generally, $\pi\in\mathbb{N}^+$ is used to avoid showing too many data points. Furthermore, any upper bound on $\pi$ can be used. However, if $\mathcal{R} = \max_{p\in\mathcal{P},a\in\mathcal{A}} r_{p,a}$,  $\forall r_{p,a} < \infty$, the performance profiles will remain the same for any $\pi \geq \mathcal{R}$. Therefore, $\mathcal{R}$ is used as the upper bound on $\pi$.

It is interesting to note that $\rho_{a}(1)$ corresponds to the percentage of problems for which algorithm $a$ has the best performance. Also, $\rho_{a}(\mathcal{R})$ represents the percentage of problems that algorithm $a$ was able to solve under the condition of the convergence test $\mathcal{C}$. Therefore, algorithms with high values of $\rho_{a}(\pi)$ are of interest.  

\subsubsection{Interpretation}
\label{sec:interp_perf_prof}

A performance profile has the values of $\pi$ on the x-axis ranging from 1 to $\mathcal{R}$, the maximum of all ratios. The proportion $\rho_a(\pi)$ is situated on the y-axis. There are three specific elements to analyze to interpret these profiles:
\begin{itemize}
    \item[-] the proportion for each algorithm $a$ at the value $\pi = 1$, \emph{i.e.} at the leftmost side of the graph. Indeed, the proportion $\rho_{a}(1)$ indicates the percentage of problems for which an algorithm $a$ is the best, based on the performance test $t_{p,a}$.
    \item[-] the proportion for each algorithm $a$ at the value $\pi = \mathcal{R}$, \emph{i.e.} at the rightmost side of the graph. This proportion indicates the percentage of problems that algorithm $a$ was able to solve within the convergence criterion $\mathcal{C}_{p,a}$.
    \item[-] how quickly the algorithm $a$ reaches a proportion of 100\%, \emph{i.e.} the line reaches the top of the graph. The value of $\pi$ for which algorithm $a$ reaches 100\% indicates its worse relative performance compared to all other algorithms. 
\end{itemize}
This thus means that the ideal performance profile starts at 100\% and finish at 100\%. It means that this algorithm $a$ is the best for the performance measure $t_{p,a}$. However, these kinds of results are rare. Therefore, it is important to compare the algorithms by looking at the three points cited above to determine which algorithm is better.

\subsubsection{Results}
\label{sec:res_perf_prof}

In our case, the set of problems $\mathcal{P}$ contains the ten models in Table~\ref{tab:models} and the set of optimization algorithms $\mathcal{A}$ include the fifteen algorithms in Table~\ref{tab:algorithms}. The convergence test $\mathcal{C}_{p,a}$ tells us whether the algorithm was able to converge with the required precision $\varepsilon$, in less than 1000 epochs. We selected two performance measures: the execution time and the number of epochs. 

\begin{figure}[H]
    \centering
    \includegraphics[width=\textwidth]{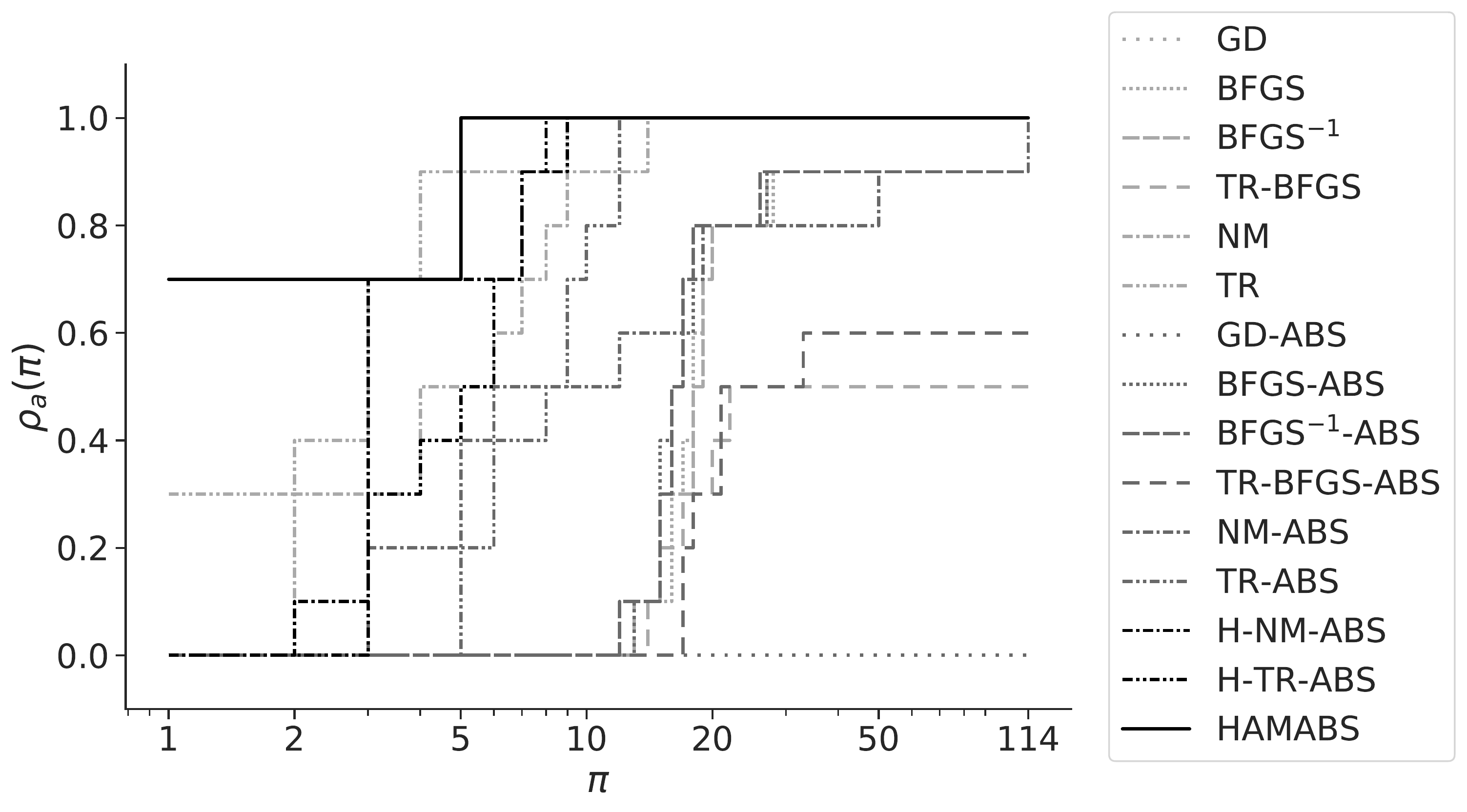}
    \caption{Performance profile on the execution time for the all models in Table~\ref{tab:models} and all algorithms in Table~\ref{tab:algorithms}.}
    \label{fig:perf_time}
\end{figure}

Figure~\ref{fig:perf_time} shows the results for the execution time. While analyzing the lines in Figure~\ref{fig:perf_time}, it is recommended to have a look at the reported execution times in Tables~\ref{tab:results_time_1} and \ref{tab:results_time_2}, provided in Appendix \ref{sec:appendix_perf_time}. We also provide the maximum ratio, $\mathcal{R}= 114$, for the execution time.

\begin{figure}[H]
    \centering
    \begin{subfigure}[b]{0.49\textwidth}
        \centering
        \includegraphics[width=\textwidth]{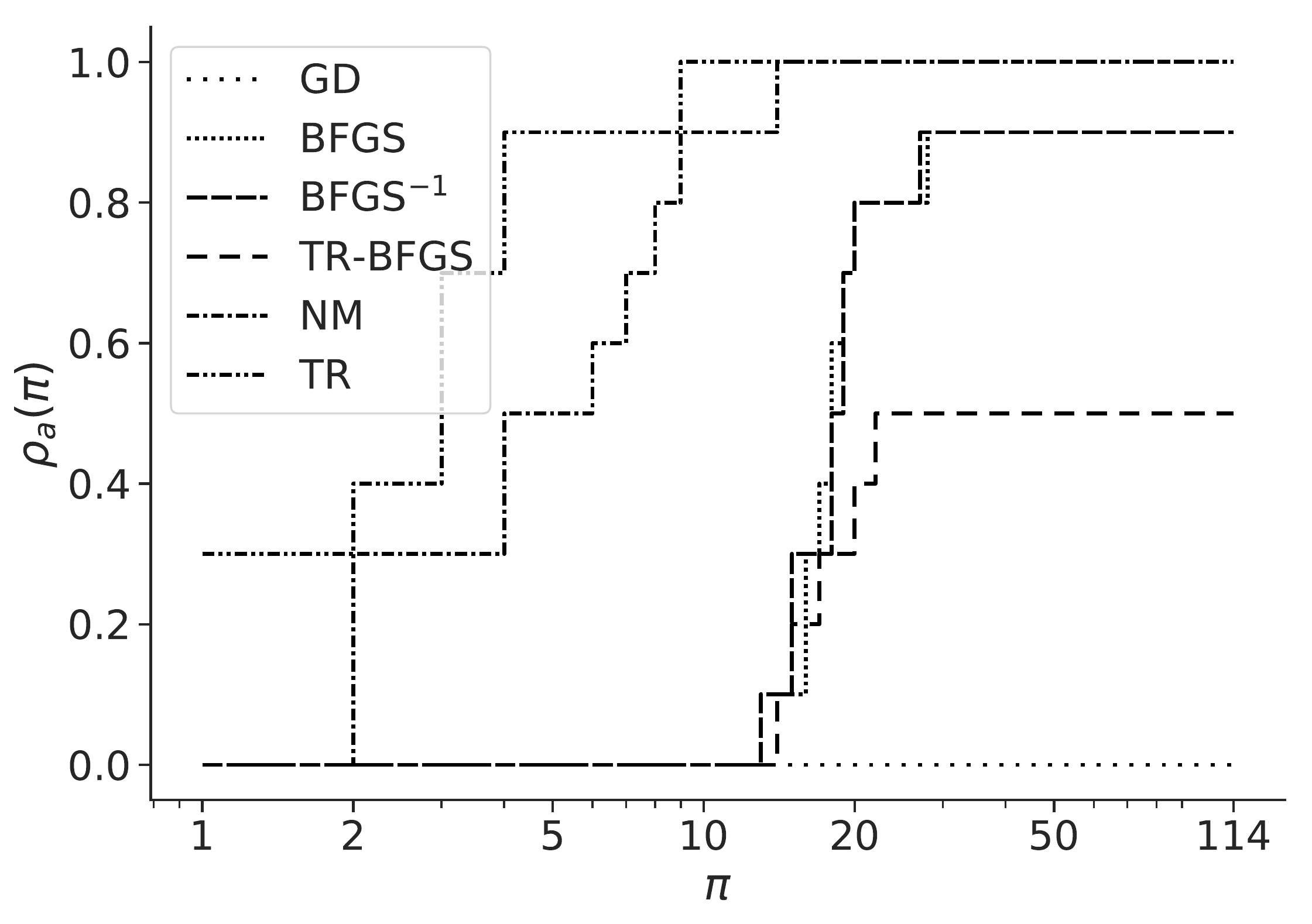}
        \caption{Standard non-stochastic algorithms}
        \label{fig:perf_profile_standard}
    \end{subfigure}
    \hfill
    \begin{subfigure}[b]{0.49\textwidth}
        \centering
        \includegraphics[width=\textwidth]{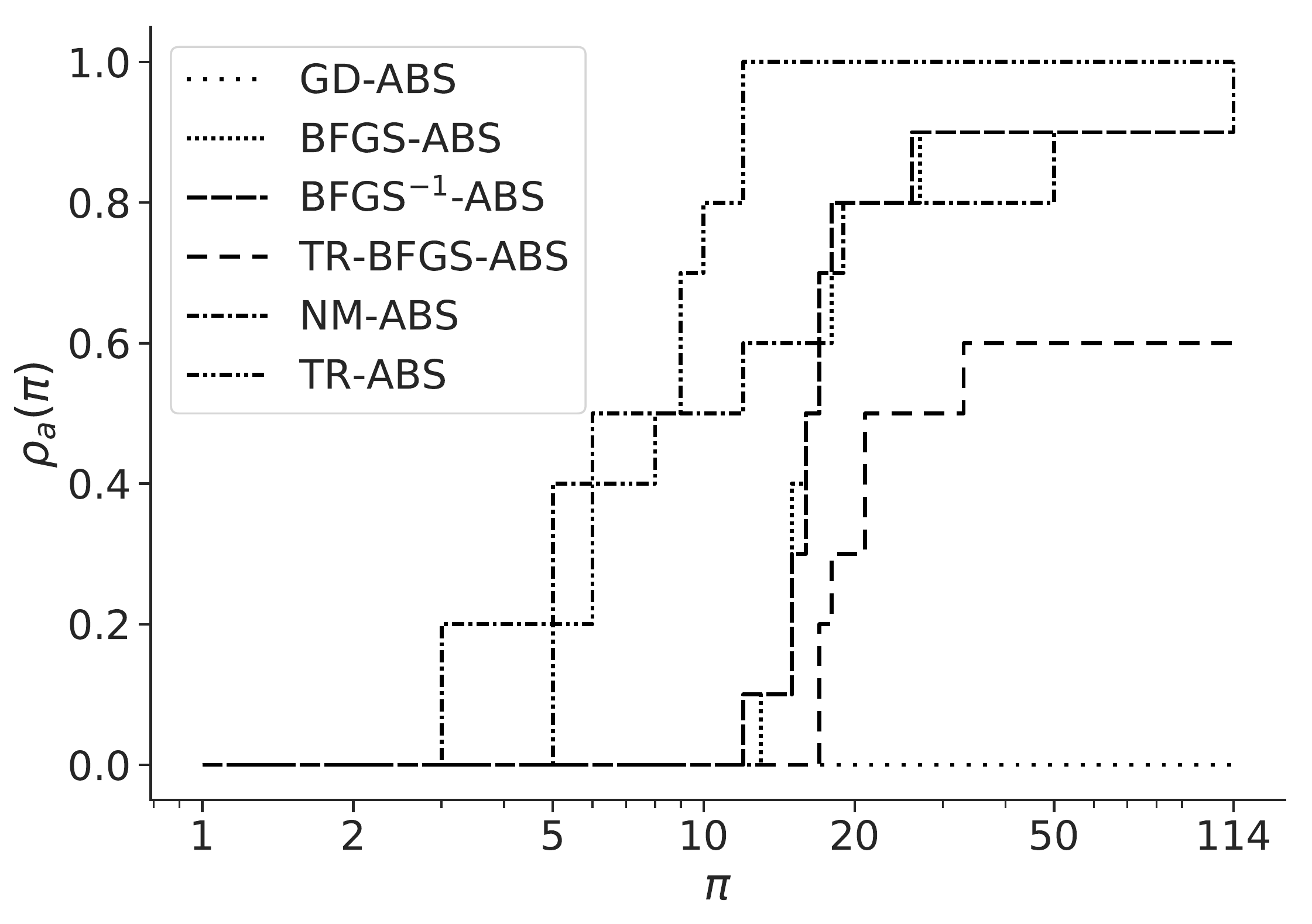}
        \caption{Stochastic \gls{amabs} algorithms}
        \label{fig:perf_profile_abs}
    \end{subfigure}
    \hfill
    \begin{subfigure}[b]{0.49\textwidth}
        \centering
        \includegraphics[width=\textwidth]{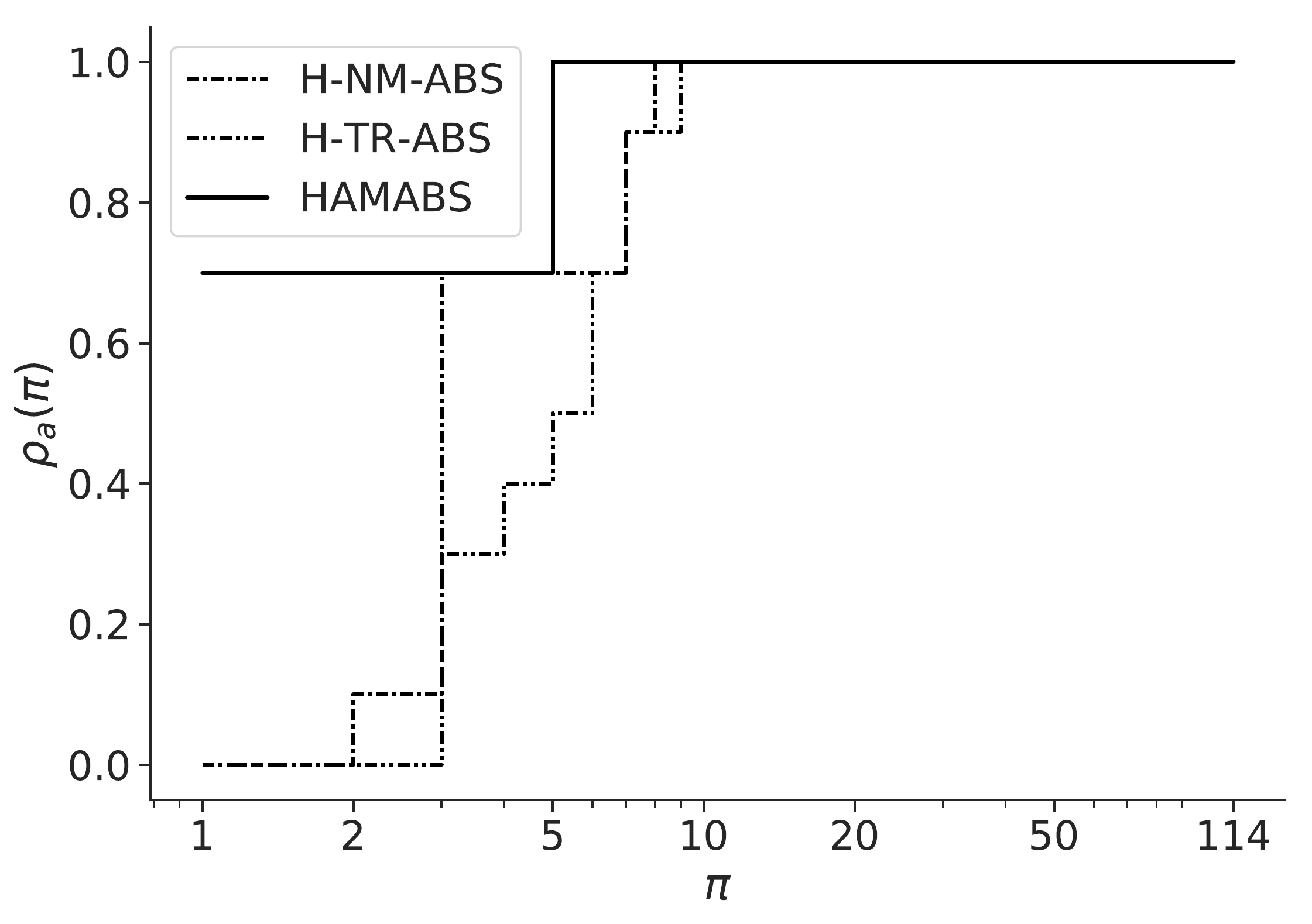}
        \caption{Hybrid stochastic algorithms}
        \label{fig:perf_profile_hybrid-abs}
    \end{subfigure}
    \hfill
    \caption{Performance profiles on the execution time for all models in Table~\ref{tab:models} splitted into different groups of algorithms.}
    \label{fig:performance_split}
\end{figure}

As an example, we analyze in detail the line corresponding to the algorithm \texttt{HAMABS} in Figure~\ref{fig:perf_time} based on the three elements cited in Section~\ref{sec:interp_perf_prof}. 
\begin{itemize}
    \item[-] it reaches a proportion of 70\% at $\pi = 1$. It, therefore, means that this algorithm is the fastest for 7 out of the 10 optimized models. Tables~\ref{tab:results_time_1} and \ref{tab:results_time_2}, provided in Appendix \ref{sec:appendix_perf_time}, report the average time and the standard deviation to optimize each model with each algorithm. As seen in these tables, the \texttt{HAMABS} is effectively the fastest algorithm on seven out of ten models.
    \item[-] it can solve all problems. Indeed, it reaches a proportion of 100\% for $\pi = \mathcal{R}$. We can also verify that it is effectively the case in Tables~\ref{tab:results_time_1} and \ref{tab:results_time_2}.
    \item[-] it reaches a proportion of 100\% at $\pi = 5$. Thus, this algorithm has, at worst, a relative performance of 5 compared to the fastest algorithm on all the models. 
\end{itemize}
Based on the analysis above and by comparing the \texttt{HAMABS} algorithms with the other algorithms, we can conclude that this algorithm is the fastest and the most robust in general. Indeed, it is the fastest on the majority of the models. Besides, the only models on which this algorithm is not the fastest are the \texttt{LPMC\_DC} models, the smallest models in terms of parameters. Also, we see that this algorithm is the fastest to reach 100\% proportion. This thus shows that it is the most robust algorithm across all models. Since there are many algorithms to analyze, we discuss them further by types of algorithms. Figure~\ref{fig:performance_split} split the lines in Figure~\ref{fig:perf_time} by the three types of algorithms showed in Section~\ref{sec:sum_algos}.

\paragraph{Standard non-stochastic algorithms}

Figure~\ref{fig:perf_profile_standard} shows the performance profile for all standard non-stochastic algorithms. We observe that these standard algorithms are struggling to optimize the models. Trust-Region and Newton's method are the fastest and the most robust methods amongst the standard ones and reach the 100\% proportion with a relative performance of up to 15 times the fastest algorithm. The two \gls{bfgs} methods are slower than the first two methods. Besides, they also fail to optimize some models. We also see that the algorithm \texttt{TR-BFGS} is struggling to optimize the models, indicating that the algorithm \texttt{H-TR-ABS} might struggle as well. The worst method is the gradient descent algorithm since it has never converged in less than 1000 epochs.

\paragraph{Stochastic algorithms}

Figure~\ref{fig:perf_profile_abs} shows the results for the algorithms using the \gls{amabs} technique. The behavior of these methods does not differ much as the slow standard methods stay slow with the \gls{amabs} method. For example, both the \texttt{NM-AMABS} and the \texttt{TR-AMABS} are the fastest and most robust algorithms. The algorithm \texttt{GD-ABS} is unable to optimize any model within the required number of epochs. Table~\ref{tab:results_time_1} and \ref{tab:results_time_2} also show that, except for the \texttt{TR-ABS}, all \gls{amabs} methods are faster than the standard ones. It thus shows that the \gls{amabs} algorithm is generally able to speed up the standard algorithms.

\paragraph{Hybrid stochastic algorithms}

Figure~\ref{fig:perf_profile_hybrid-abs} shows the results for the three algorithms using the Hybridization and the \gls{amabs} method. These three methods are the fastest algorithms on larger models. While we already discussed the case of the \texttt{HAMABS} algorithm, the other two methods are never the fastest method for any models. However, they are slightly more robust than the other algorithms. Indeed, the \texttt{H-NM-ABS} is almost as good as the \texttt{HAMABS} algorithm. However, by looking at the times, we still see quite a difference between these two algorithms. Indeed, there is generally a 20\% difference in execution time between these two algorithms. The \texttt{H-TR-ABS} seems to perform quite well. However, by looking at the times, we see that both the \texttt{TR} and the \texttt{TR-ABS} algorithms are faster. This is most likely due to the use of the Trust-Region method with the \gls{bfgs} approximation being exceptionally slow. 

\begin{figure}[H]
    \centering
    \includegraphics[width=\textwidth]{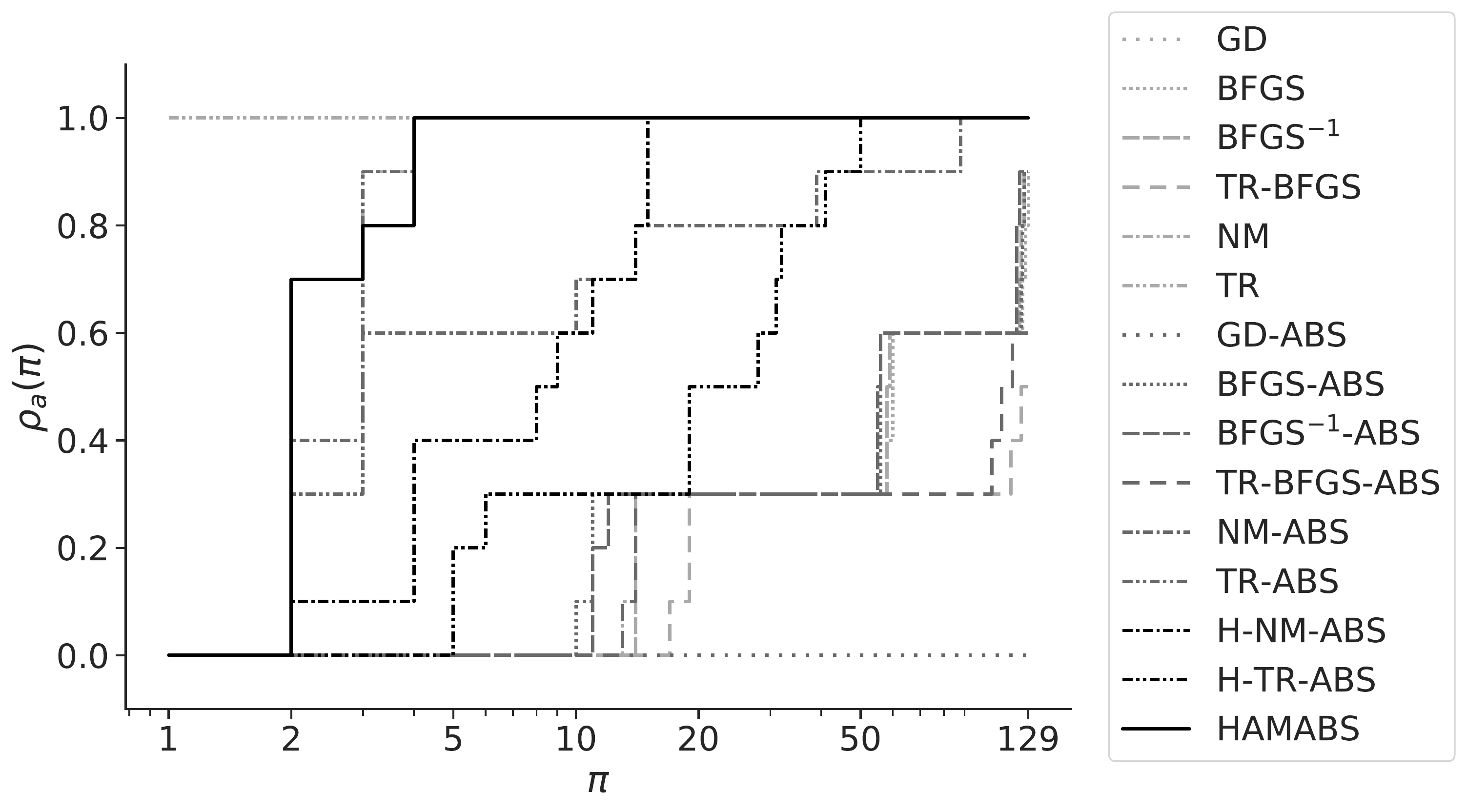}
    \caption{Performance profile on the epochs for the all models in Table~\ref{tab:models} and all algorithms in Table~\ref{tab:algorithms}.}
    \label{fig:perf_epochs}
\end{figure}

\begin{figure}[H]
    \centering
    \begin{subfigure}[b]{0.49\textwidth}
        \centering
        \includegraphics[width=\textwidth]{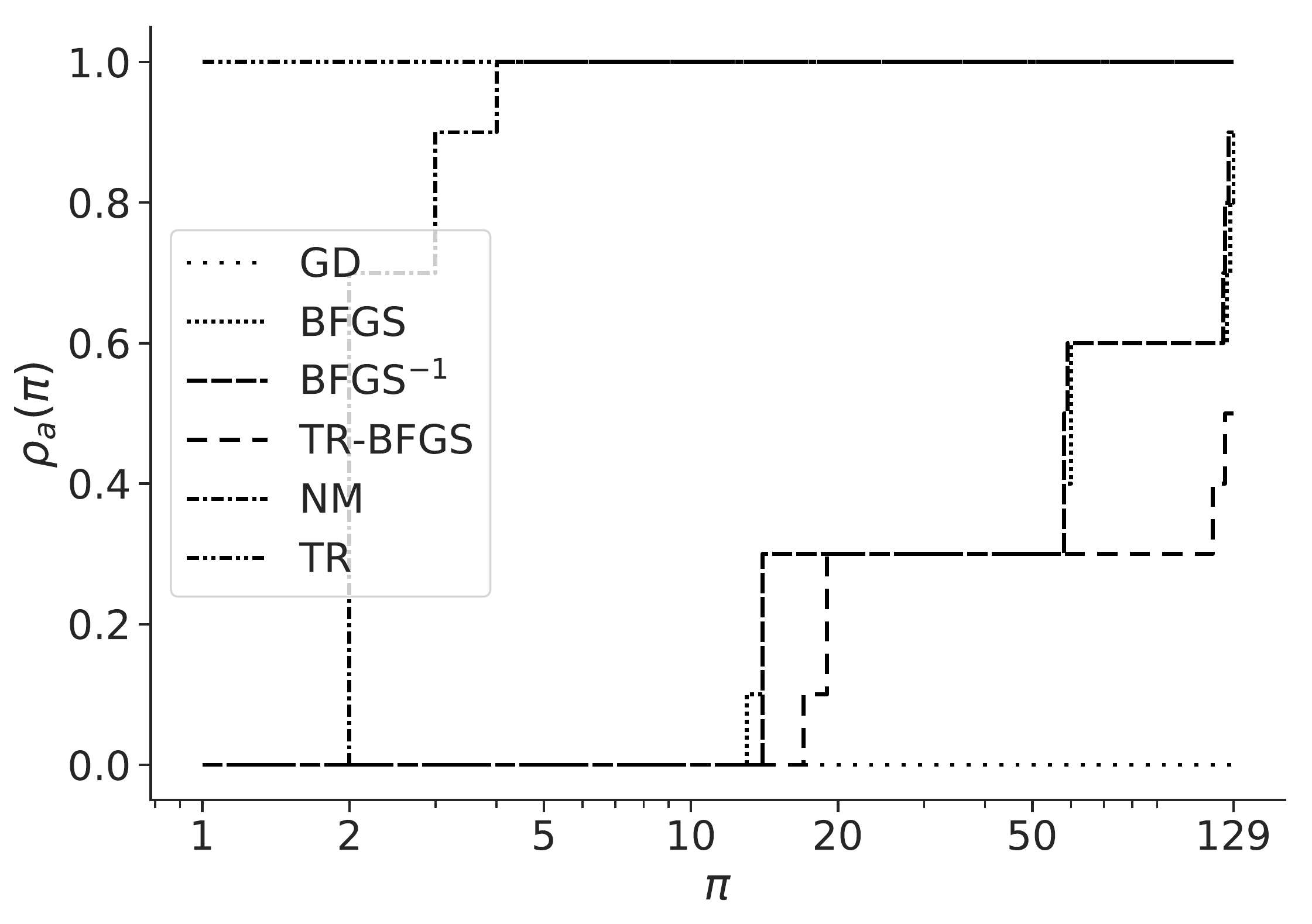}
        \caption{Standard non-stochastic algorithms}
        \label{fig:perf_profile_standard_epochs}
    \end{subfigure}
    \hfill
    \begin{subfigure}[b]{0.49\textwidth}
        \centering
        \includegraphics[width=\textwidth]{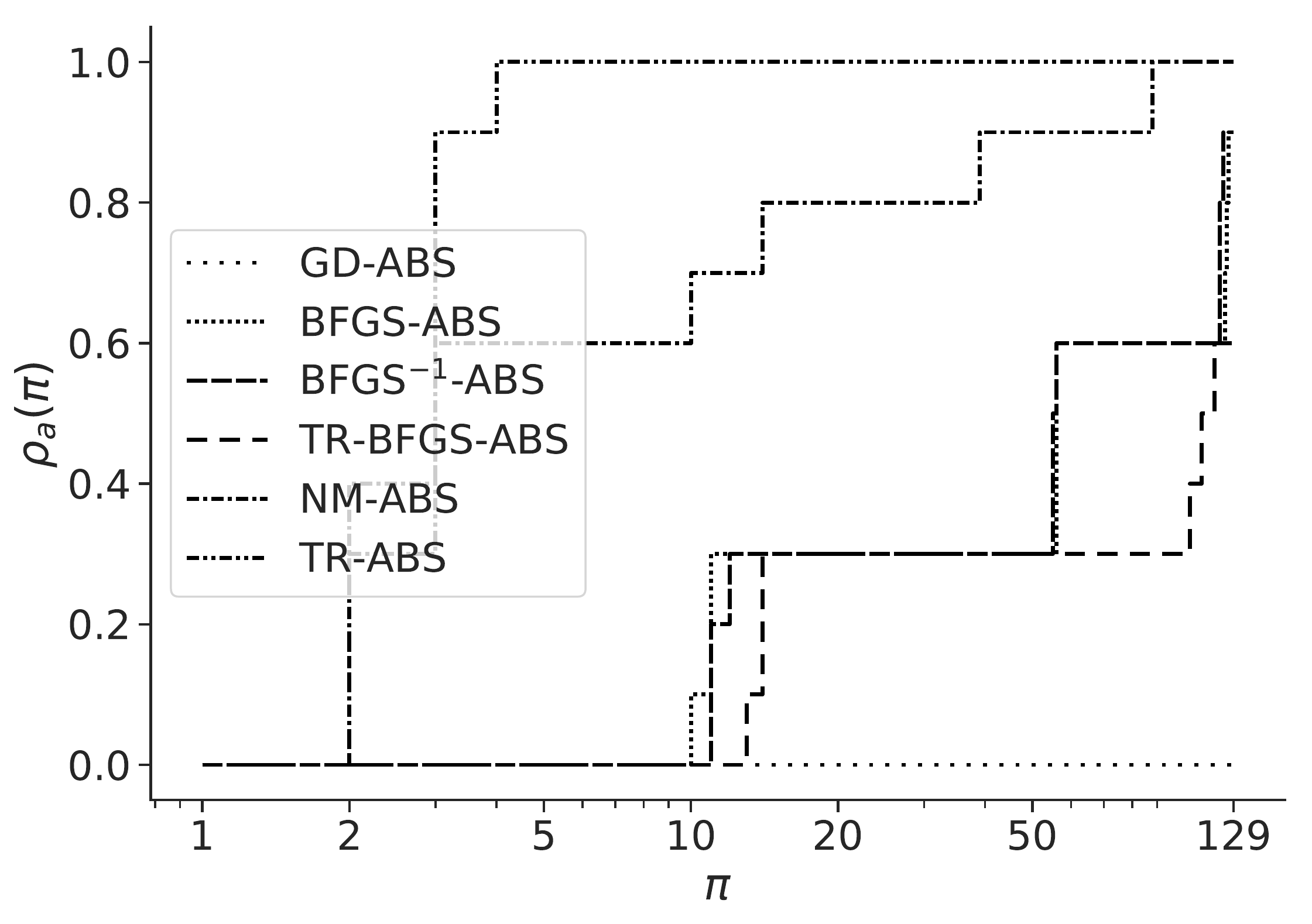}
        \caption{Stochastic \gls{amabs} algorithms}
        \label{fig:perf_profile_abs_epochs}
    \end{subfigure}
    \hfill
    \begin{subfigure}[b]{0.49\textwidth}
        \centering
        \includegraphics[width=\textwidth]{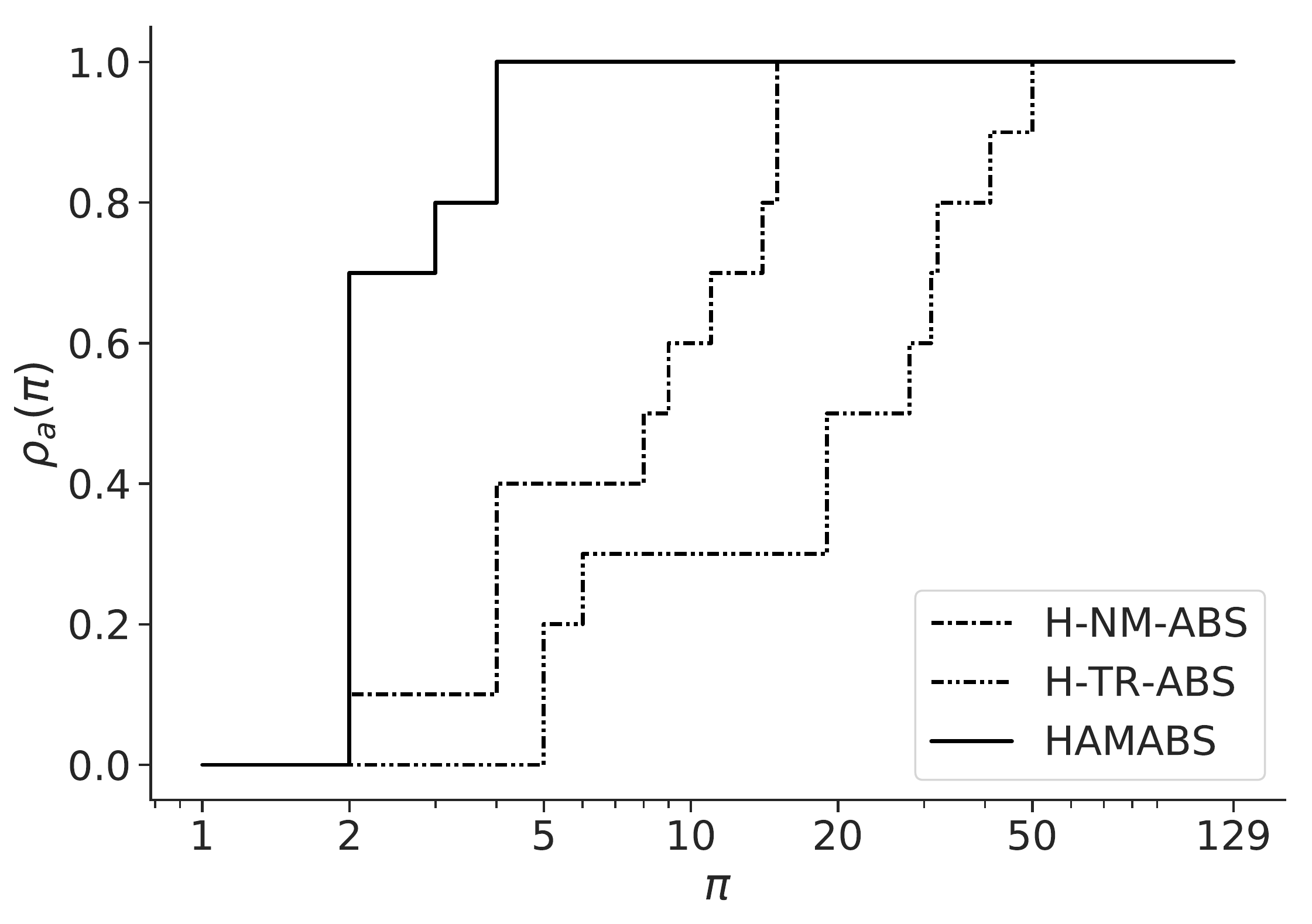}
        \caption{Hybrid stochastic algorithms}
        \label{fig:perf_profile_hybrid-abs_epochs}
    \end{subfigure}
    \hfill
    \caption{Performance profiles on the epochs for all models in Table~\ref{tab:models} splitted into different groups of algorithms.}
    \label{fig:performance_split_epochs}
\end{figure}

As seen in Figures~\ref{fig:perf_time} and \ref{fig:performance_split}, the \texttt{HAMABS} algorithm is the fastest to optimize most of the models. While the execution time is important, the number of epochs used throughout the optimization process can also bring some important conclusions. An optimization algorithm performs one epoch as soon as it has seen all the data in the dataset once. Thus, the time to optimize a model is often correlated to the number of epochs it takes. Furthermore, this performance measure is independent of the type and the load of the computer.  Figure~\ref{fig:perf_epochs} shows the performance profile on the epochs for all algorithms, and Figure~\ref{fig:performance_split_epochs} splits the previous figure based on the different groups of algorithms. We can see in Figure~\ref{fig:perf_epochs} that the \texttt{HAMABS} algorithm is more robust than most of the other algorithms, except for the \texttt{TR} algorithm. 

If we compare the different algorithms, we see that second-order methods tend to use fewer epochs to achieve convergence. Indeed, we can see in both Figure~\ref{fig:perf_profile_standard_epochs} and Figure~\ref{fig:perf_profile_abs_epochs} that the methods based on Newton's method (\texttt{NM}/\texttt{NM-ABS}) and Trust-Region method (\texttt{TR}/\texttt{TR-ABS}) are more robust than the other algorithms. This is expected because these methods use the information on the curvature. They thus require fewer steps to complete the estimation process. We also see that second-order methods using the full size dataset, \texttt{NM} and \texttt{TR}, are using less epochs than the stochastic methods, \texttt{NM-ABS} and \texttt{TR-ABS}. This is also an expected behavior since the stochastic algorithms use less information per step. They thus need to perform many more steps, often leading to more epochs, to gain the same knowledge. On the other hand, they spend less time on each step, leading to a consequent speedup. It is interesting to note that the \texttt{HAMABS} algorithm is amongst the algorithms using the least number of epochs. Indeed, it reaches a proportion of 100\% with a relative performance of 4. It, therefore, explains why this algorithm is that fast compared to the \gls{amabs} algorithms. Also, we see that the \texttt{H-NM-ABS} tends to use more epochs than the \texttt{HAMABS} algorithm. This could explain why \texttt{HAMABS} is the fastest algorithm. Tables~\ref{tab:results_epochs_1} and \ref{tab:results_epochs_2} in Appendix \ref{sec:appendix_perf_epochs} show the average number of epochs with the standard deviation used by each algorithm to optimize each model.

\subsection{Comparison with State-of-the-Art Software}
\label{sec:comp_software}

We now compare the performance of our best algorithm, the \texttt{HAMABS} algorithm, to Pandas Biogeme~\citep{bierlaire_pandasbiogeme:_2018}, a state-of-the-art choice modeling software. Biogeme is using the Python package Scipy to optimize the models. The default algorithm for the minimization in this package is the \gls{bfgs}$^{-1}$. Table~\ref{tab:comp_biogeme_time} reports the average time to optimize each model for both Biogeme and the \texttt{HAMABS} algorithm. Besides, the last column shows the speedup gained by using the \texttt{HAMABS} algorithm instead of Scipy. 

\begin{table}[H]
  \centering
  \renewcommand\arraystretch{1.1}
  \begin{tabular}{l||c|c|c}
    \multirow{2}{*}{Models} & \multicolumn{2}{c|}{Time [s]} & \multirow{2}{*}{Speedup} \\ \cline{2-3}
    & \texttt{HAMABS} & Biogeme/Scipy & \\ \hline \hline
    \texttt{LPMC\_DC\_S} & $1.86 \pm 0.12$ & $1.62 \pm 0.01$ & $\times 1.15$ \\
    \texttt{LPMC\_DC\_M} & $3.11 \pm 0.20$ & $2.79 \pm 0.02$ & $\times 1.11$ \\
    \texttt{LPMC\_DC\_L} & $4.59 \pm 0.32$ & $4.07 \pm 0.04$ & $\times 1.13$ \\ \hline
    \texttt{LPMC\_RR\_S} & $11.98 \pm 1.23$ & $65.17 \pm 0.09$ & $\div 5.44$ \\
    \texttt{LPMC\_RR\_M} & $18.46 \pm 1.06$ & $127.67 \pm 0.30$ & $\div 6.91$ \\
    \texttt{LPMC\_RR\_L} & $18.14 \pm 1.06$ & $177.09 \pm 0.29$ & $\div 9.76$ \\ \hline
    \texttt{LPMC\_Full\_S} & $257.02 \pm 42.85$ & $1462.51 \pm 14.41$ & $\div 5.69$ \\
    \texttt{LPMC\_Full\_M} & $405.43 \pm 43.77$ & $2480.06 \pm 18.27$ & $\div 6.12$ \\
    \texttt{LPMC\_Full\_L} & $486.31 \pm 63.38$ & $4758.28 \pm 45.22$ & $\div 9.78$ \\ \hline \hline
    \texttt{MTMC} & $1243.95 \pm 56.21$ & $28008.10 \pm 528.33$ & $\div 22.52$ \\
  \end{tabular}
  \caption{Comparison of the optimization time for all models in Table~\ref{tab:models} between the \texttt{HAMABS} algorithm and Biogeme. The time are reported in seconds. The speedup corresponds to a ratio between the two compared values.}
  \label{tab:comp_biogeme_time}
\end{table}

Results presented in Table~\ref{tab:comp_biogeme_time} show that the algorithm \texttt{HAMABS} is generally faster than the Scipy package. On the models \texttt{\gls{lpmc}\_DC}, that have few parameters, the \texttt{HAMABS} algorithm is slower with a ratio around 1.15. However, the \texttt{HAMABS} becomes faster than the Scipy package on the models \texttt{\gls{lpmc}\_RR} and \texttt{\gls{lpmc}\_Full} that include more parameters. This implies that a model that previously took minutes, or hours to converge, is now optimized in only a few seconds or minutes, respectively. The most important gain is on the optimization of the largest model, the \texttt{\gls{mtmc}} model, with a speedup ratio exceeding 22. While Biogeme took around seven and a half hours to converge, our \texttt{HAMABS} algorithm is converging in less than 20 minutes. 

Table~\ref{tab:comp_biogeme_epochs} compares the two algorithms based on the number of epochs used for optimization. Results show that the computational gains achieved by our \texttt{HAMABS} algorithm are even more important in terms of the number of epochs. For Biogeme and the Scipy package, the number of epochs is directly correlated to the number of parameters. As a result, the number of epochs highly depends on the size of the model. For example, the \texttt{\gls{mtmc}} models uses hundred times more epochs to be optimized than the \texttt{\gls{lpmc}\_DC} models.  For our \texttt{HAMABS} algorithm, on the other hand, the number of epochs used is more stable across the different models. Indeed, between the smallest and the largest models, the average number of epochs only doubles. On the \texttt{\gls{mtmc}} model, the speedup ratio in number of epochs between Biogeme and our \texttt{HAMABS} algorithm exceeds 600. The reason behind these discrepancies in the number of epochs is the stopping criterion. The Scipy package is using the standard, yet incorrect, gradient value as the stopping criterion. If the objective function and its derivatives are not correctly normalized, this criterion can either stop the algorithm too early or too late. Therefore, the use of an appropriate stopping criterion makes an important difference in terms of epochs while keeping sufficient precision, as seen in the last column of Table~\ref{tab:comp_biogeme_epochs}. Indeed, we see that the relative difference in log likelihood between the two methods is never larger than $2\times 10^{-4}$ \%. 

\begin{table}[H]
  \centering
  \renewcommand\arraystretch{1.1}
  \begin{tabular}{l||c|c|c|c}
    \multirow{2}{*}{Models} & \multicolumn{2}{c|}{Epochs} & \multirow{2}{*}{Speedup} & \multirow{2}{*}{$\Delta \mathcal{L}$ [\%]} \\ \cline{2-3}
    & \texttt{HAMABS} & Biogeme/Scipy & & \\ \hline \hline
    \texttt{LPMC\_DC\_S} & $13.79 \pm 1.70$ & $122$ & $\div 8.85$ & $3.20\times10^{-6}$ \\
    \texttt{LPMC\_DC\_M} & $13.50 \pm 2.05$ & $114$ & $\div 8.44$ & $3.13\times10^{-6}$ \\
    \texttt{LPMC\_DC\_L} & $12.57 \pm 1.93$ & $123$ & $\div 9.78$ & $5.43\times10^{-6}$ \\ \hline
    \texttt{LPMC\_RR\_S} & $15.31 \pm 2.53$ & $787$ & $\div 51.40$ & $5.00\times10^{-5}$ \\
    \texttt{LPMC\_RR\_M} & $13.95 \pm 1.67$ & $809$ & $\div 57.97$ & $3.28\times10^{-5}$ \\
    \texttt{LPMC\_RR\_L} & $9.55 \pm 0.56$ & $772$ & $\div 80.84$ & $2.87\times10^{-4}$ \\ \hline
    \texttt{LPMC\_Full\_S} & $24.98 \pm 2.16$ & $1786$ & $\div 71.50$ & $1.45\times10^{-4}$ \\
    \texttt{LPMC\_Full\_M} & $20.42 \pm 1.84$ & $1531$ & $\div 74.96$ & $2.13\times10^{-4}$ \\
    \texttt{LPMC\_Full\_L} & $21.36 \pm 2.05$ & $1996$ & $\div 93.44$ & $8.02\times10^{-5}$ \\ \hline \hline
    \texttt{MTMC} & $18.63 \pm 1.58$ & $11920$ & $\div 639.88$ & $2.01\times10^{-4}$ \\
  \end{tabular}
  \caption{Comparison of the epochs used in the optimization process for all models in Table~\ref{tab:models} between the \texttt{HAMABS} algorithm and Biogeme. The time are reported in seconds. The speedup corresponds to a ratio between the two compared values. The last column corresponds to the relative difference in percentage between the final value of the log likelihood returned by Biogeme and by the \texttt{HAMABS} algorithm.}
  \label{tab:comp_biogeme_epochs}
\end{table}

The similar ratios between the models \texttt{\gls{lpmc}\_RR} and \texttt{\gls{lpmc}\_Full}, can be due to the added complexity on the \texttt{\gls{lpmc}\_Full} models. Indeed, in these models, multiple parameters are computed on small populations. This leads to an increase in complexity, and the stochasticity might, therefore, not be that helpful. The algorithm has to perform more steps at the full size to find the parameter values on the small groups. We thus lose some time at the end of the optimization process compared to \texttt{\gls{lpmc}\_RR} models.

In definitive, our results showed that the \texttt{HAMABS} algorithm is not only the fastest among the 15 algorithms in Table~\ref{tab:algorithms}, but that it is also much faster than the current implementation of the state-of-the-art choice modeling software Biogeme.

\subsection{Sensitivity Analysis}
\label{sec:parameters_algo}
We want now to test the sensitivity of the \texttt{HAMABS} algorithm's parameters to make sure that validate the choice of the default parameters given in Algorithm~\ref{algo:hybrid}. We selected the three following models to perform the sensitivity study: \texttt{\gls{lpmc}\_DC\_L}, \texttt{\gls{lpmc}\_RR\_L}, and \texttt{\gls{lpmc}\_Full\_L}.
The sensitivity analysis was performed on the estimation time of these models by the \texttt{HAMABS} algorithm. Each model was trained 20 times for all the values present in Table~\ref{tab:param_sens}. 

\begin{table}[H]
  \centering
  \renewcommand\arraystretch{1.1}
  \begin{tabular}{l||c|c}
    Parameter & Default & Test values \\ \hline\hline
    $W$ & 10 & [1, 2, 3, $\cdots$, 18, 19, 20] \\
    $\Delta$ & 1\% & [0.1, 0.2, 0.5, 1, 2, 5, 10, 20, 50, 100] \\
    $C$ & 2 & [1, 2, 3, $\cdots$, 13, 14, 15] \\
    $\tau$ & 2 & [1.1, 1.2, 1.3, 1.4, 1.5, 2, 2.5, 3, 3.5, 4, 4.5, 5, 6, 7, 8, 9, 10] \\
    $\Delta_H$ & 30\% & [0, 5, 10, $\cdots$, 90, 95, 100] \\
    $\varepsilon$ & $10^{-6}$ & [$10^{-9}$, $10^{-8}$, $\cdots$, $10^{-1}$, $10^{0}$]
  \end{tabular}
  \caption{Parameter values of the \texttt{HAMABS} algorithm used for the sensitivity analysis. We refer the reader to Algorithm~\ref{algo:hybrid} for a detailed explanation of the parameters.}
  \label{tab:param_sens}
\end{table}

All results are reported using graphs for which the relative performance is indicated on the vertical axis. All performances are normalized to the execution time obtained with the default parameter values (base value of 1) to ease the comparison of results across the models. 

\begin{figure}[H]
    \centering
    \begin{subfigure}[b]{0.49\textwidth}
        \centering
        \includegraphics[width=\textwidth]{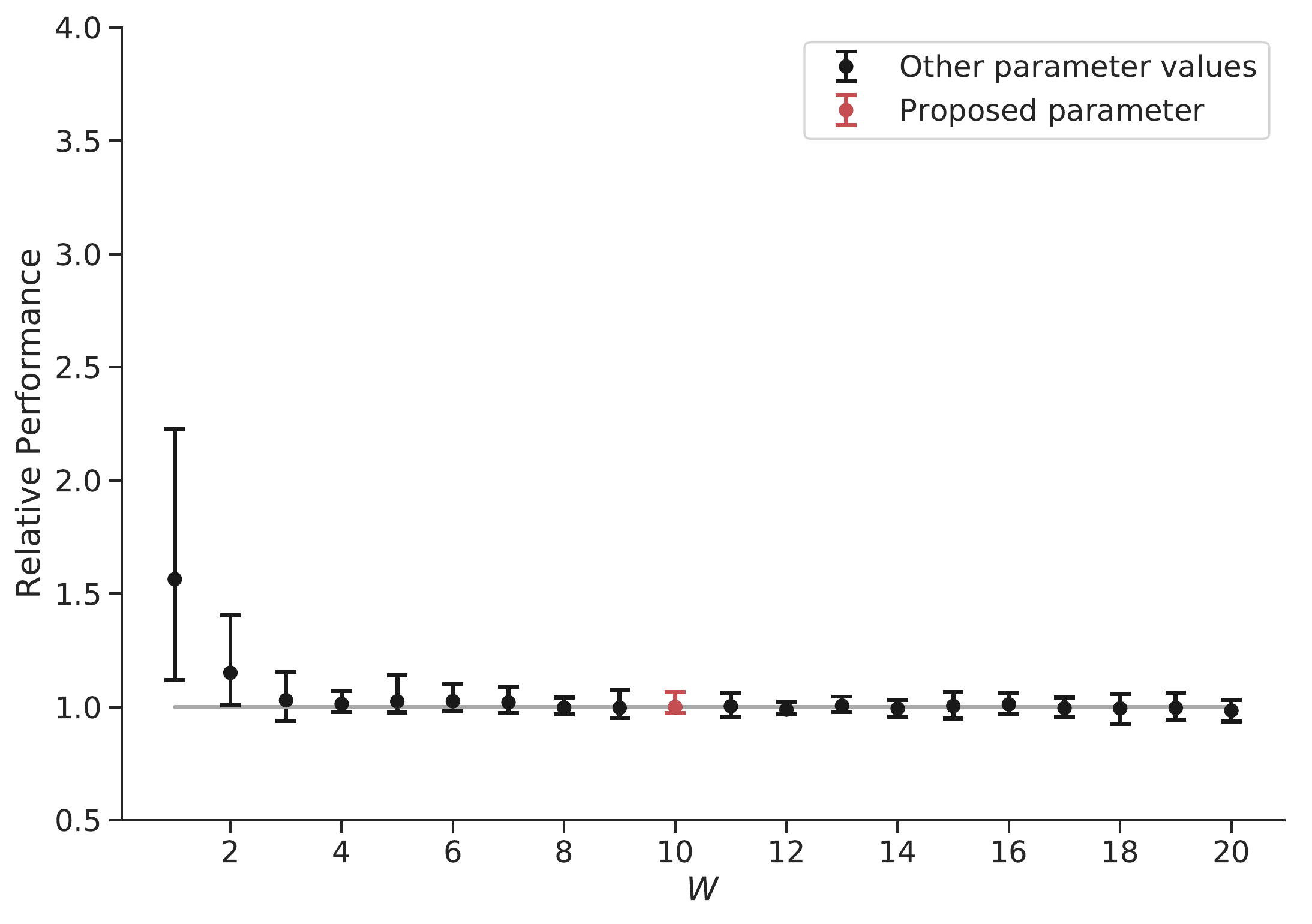}
        \caption{\texttt{\gls{lpmc}\_DC\_L}}
        \label{fig:window_DC}
    \end{subfigure}
    \hfill
    \begin{subfigure}[b]{0.49\textwidth}
        \centering
        \includegraphics[width=\textwidth]{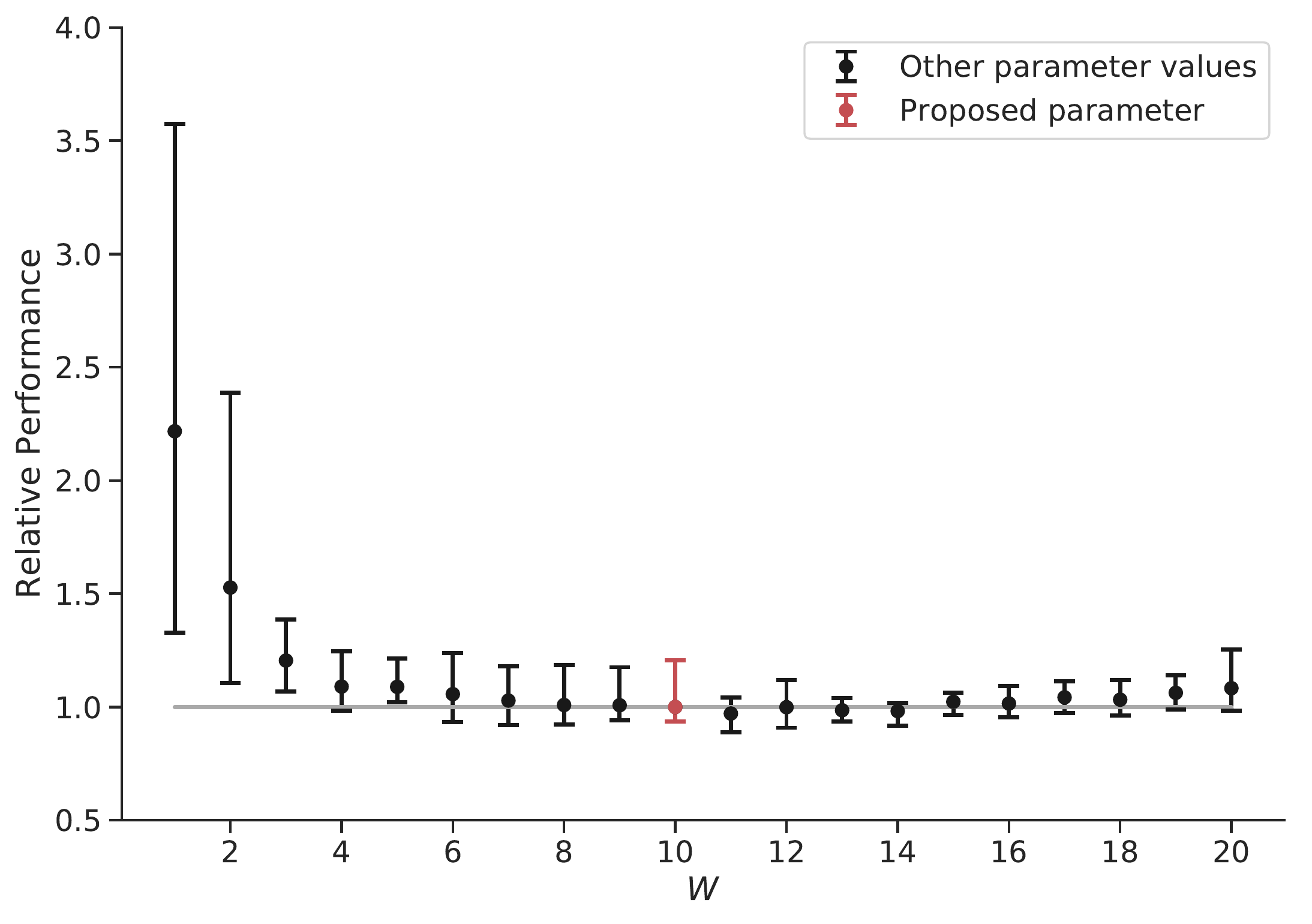}
        \caption{\texttt{\gls{lpmc}\_RR\_L}}
        \label{fig:window_RR}
    \end{subfigure}
    \hfill
    \begin{subfigure}[b]{0.49\textwidth}
        \centering
        \includegraphics[width=\textwidth]{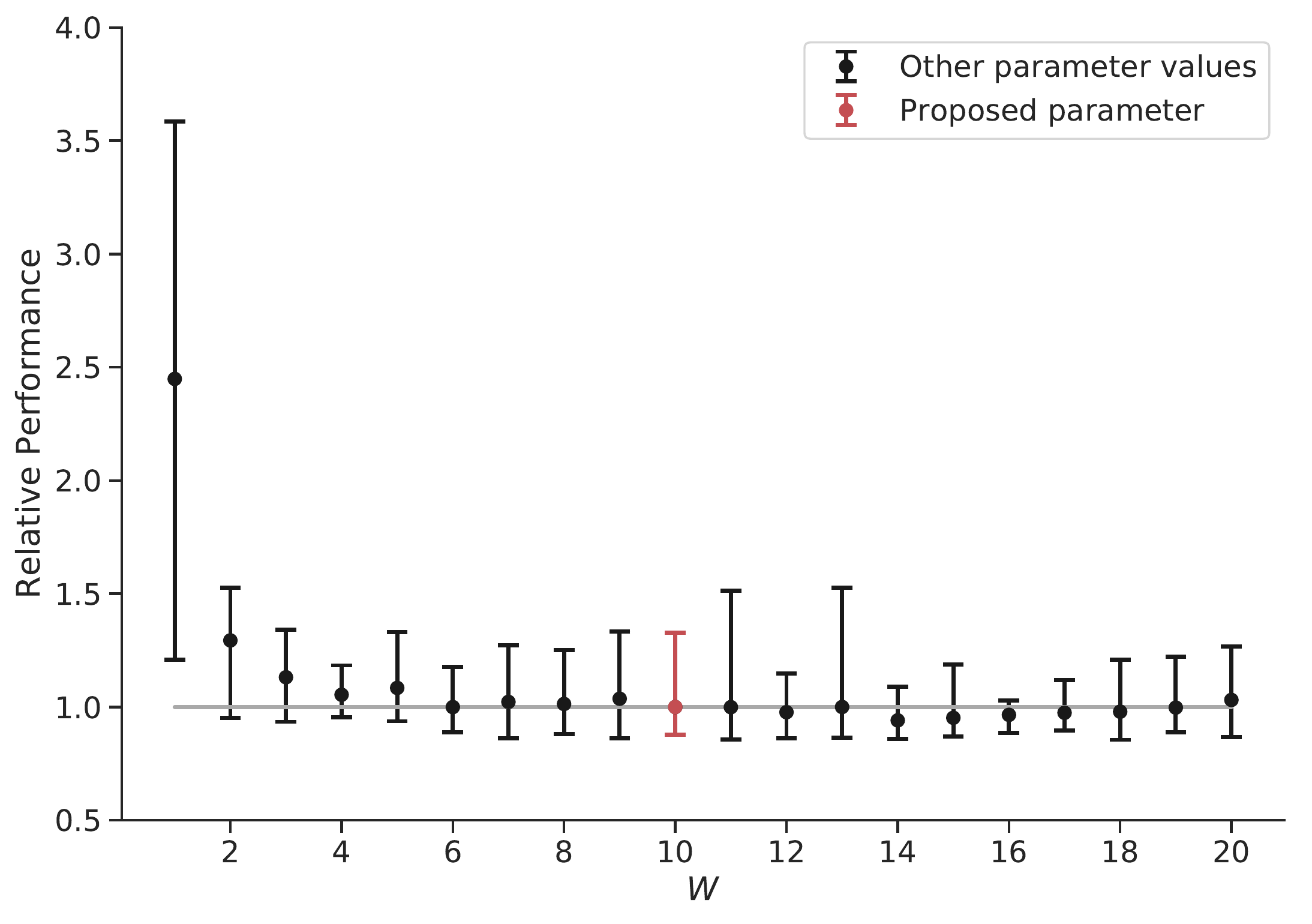}
        \caption{\texttt{\gls{lpmc}\_Full\_L}}
        \label{fig:window_Full}
    \end{subfigure}
    \hfill
    \caption{Sensitivity analysis for the parameter $W$ for the three \gls{lpmc} models using the large dataset. The black error bars correspond to the tested values and the red to the proposed value. The gray line correspond to the benchmark with the proposed parameters for the relative performance.}
    \label{fig:effect_window}
\end{figure}

Figure~\ref{fig:effect_window} shows the analysis for the parameter $W$; the size of the window. The results are similar across all models. Indeed, small values of $W$ lead to an increase in the execution time. It can be explained by the fact that if the window is too small, the average is then noisy. Therefore, the computation of the improvement of the log likelihood is not precise. It thus increases the batch size too soon or too late. This thus leads to an increase in the total execution time. Large values for $W$ do not affect the optimization process as much as small values. Indeed, we see a small increase in the execution time in Figure~\ref{fig:window_RR} when $W$ is large. Indeed, if we use large values of $W$, the average is less influenced by the new data points. Thus, the \texttt{AMABS} reaction is slower. Therefore, a good value for $W$ has to be in the middle. We thus propose to use $W = 10$. 

\begin{figure}[H]
    \centering
    \begin{subfigure}[b]{0.49\textwidth}
        \centering
        \includegraphics[width=\textwidth]{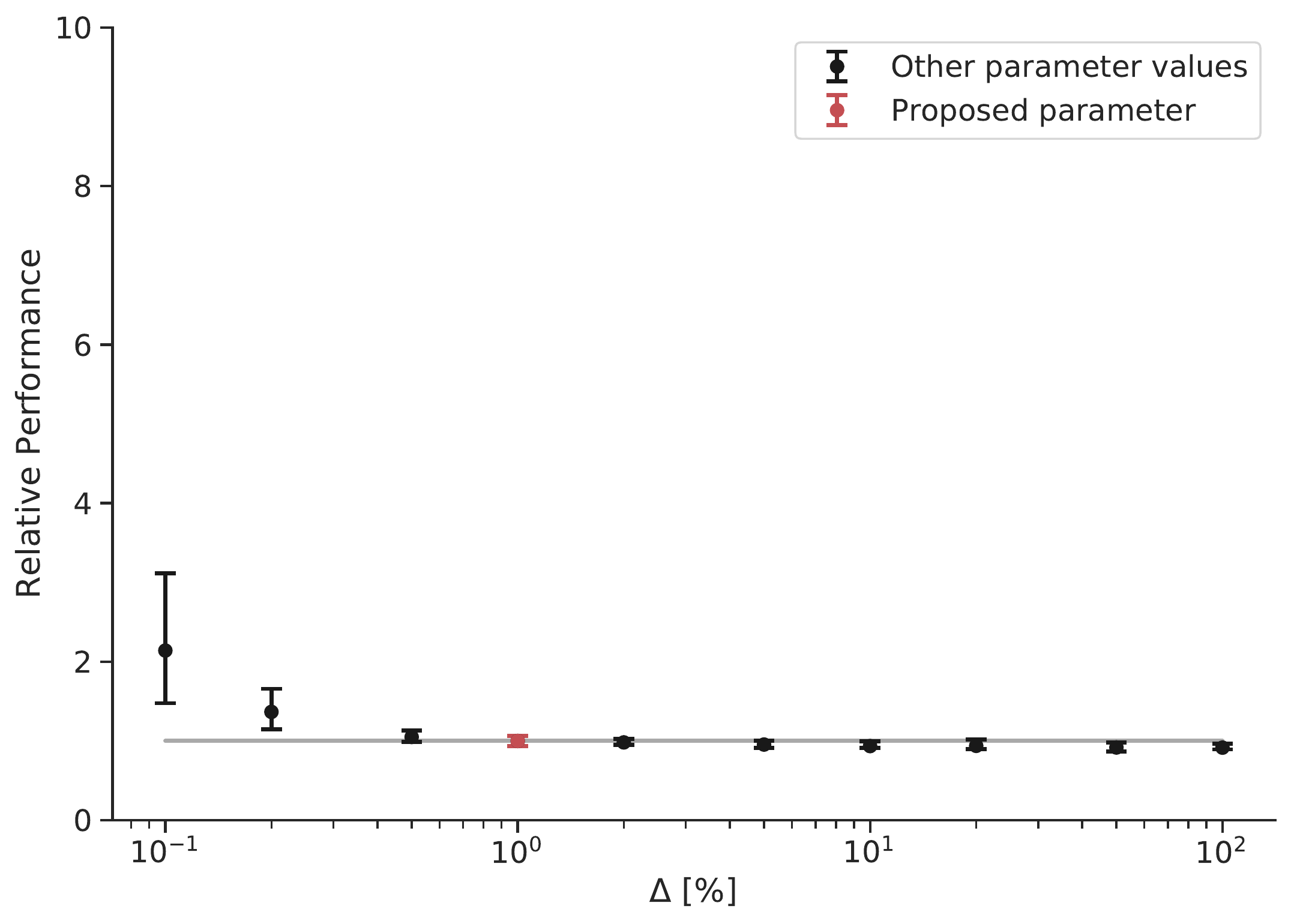}
        \caption{\texttt{\gls{lpmc}\_DC\_L}}
        \label{fig:thresh_DC}
    \end{subfigure}
    \hfill
    \begin{subfigure}[b]{0.49\textwidth}
        \centering
        \includegraphics[width=\textwidth]{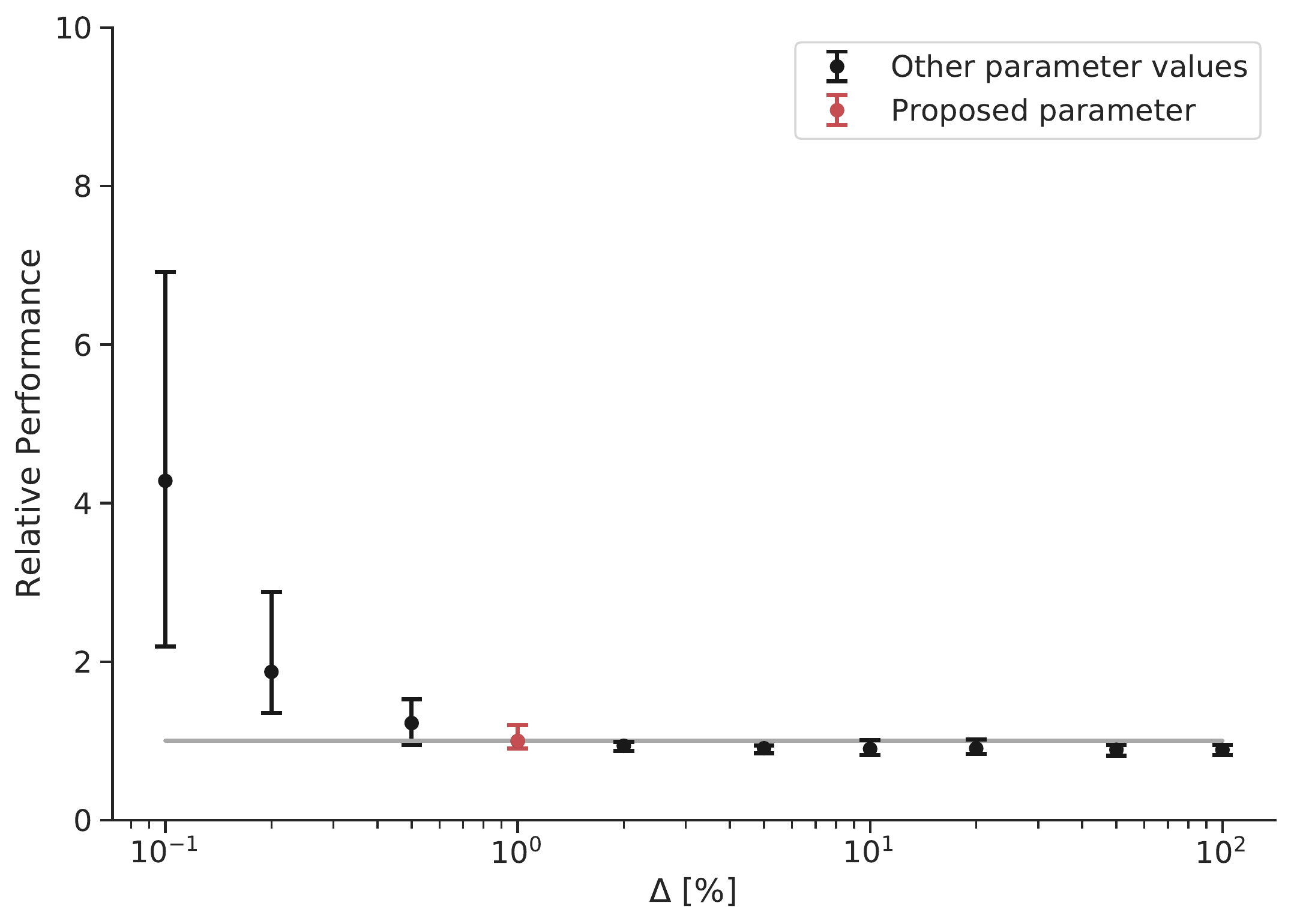}
        \caption{\texttt{\gls{lpmc}\_RR\_L}}
        \label{fig:thresh_RR}
    \end{subfigure}
    \hfill
    \begin{subfigure}[b]{0.49\textwidth}
        \centering
        \includegraphics[width=\textwidth]{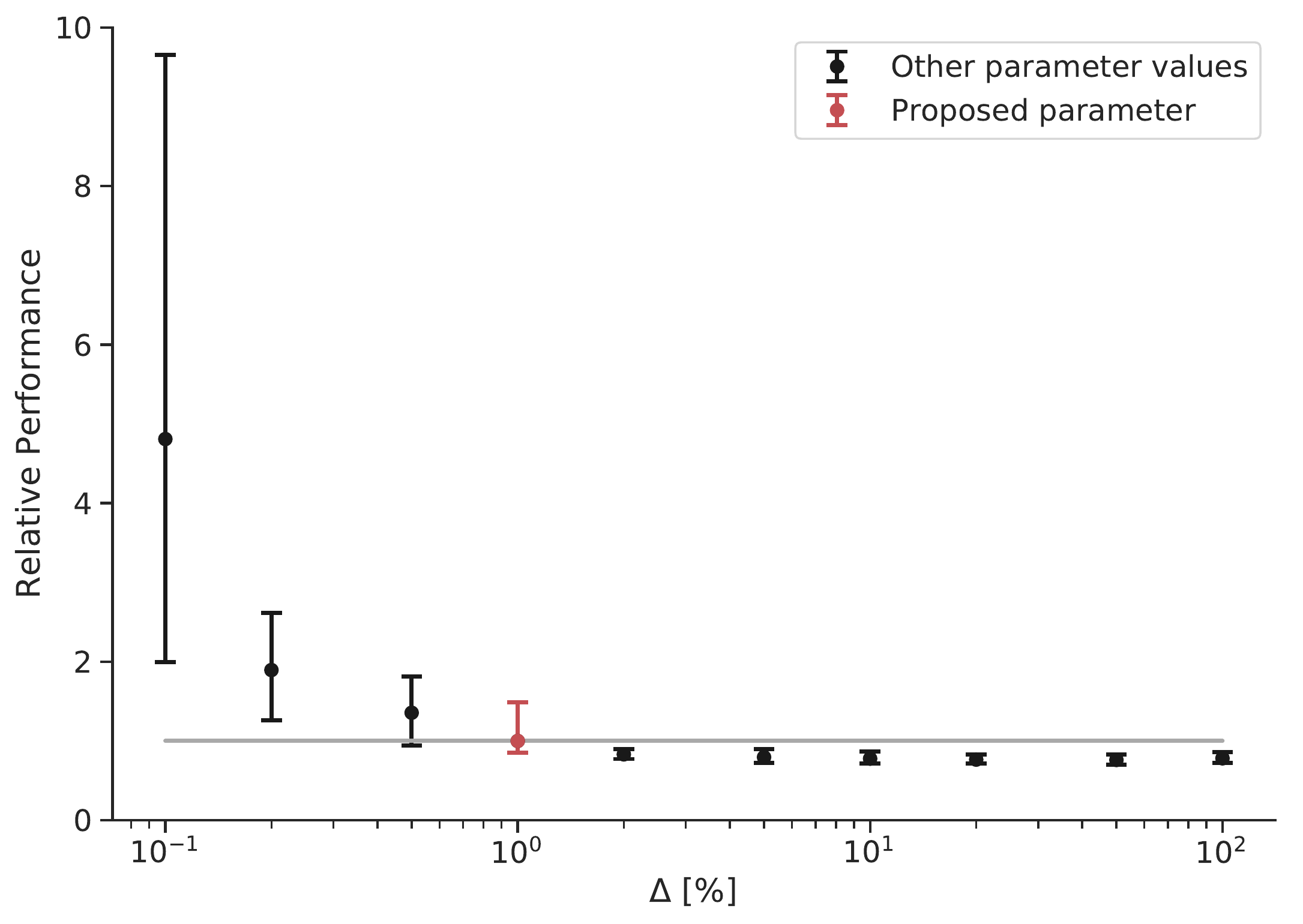}
        \caption{\texttt{\gls{lpmc}\_Full\_L}}
        \label{fig:thresh_Full}
    \end{subfigure}
    \hfill
    \caption{Sensitivity analysis for the parameter $\Delta$ for the three \gls{lpmc} models using the large dataset. The black error bars correspond to the tested values and the red to the proposed value. The gray line correspond to the benchmark with the proposed parameters for the relative performance.}
    \label{fig:effect_thresh}
\end{figure}

Figure~\ref{fig:effect_thresh} shows the analysis for the parameter $\Delta$; the threshold for successful iterations. Small values of this parameter lead to more substantial execution time. Indeed, the \texttt{HAMABS} algorithm is delaying the update of the batch size for too long. While larger values seem to perform slightly better, it also increases the probability of switching too soon for larger models. However, the value 1\% corresponds to the moment where the curve is becoming a plateau. We thus propose to use this value for the parameter $\Delta$.

\begin{figure}[H]
    \centering
    \begin{subfigure}[b]{0.49\textwidth}
        \centering
        \includegraphics[width=\textwidth]{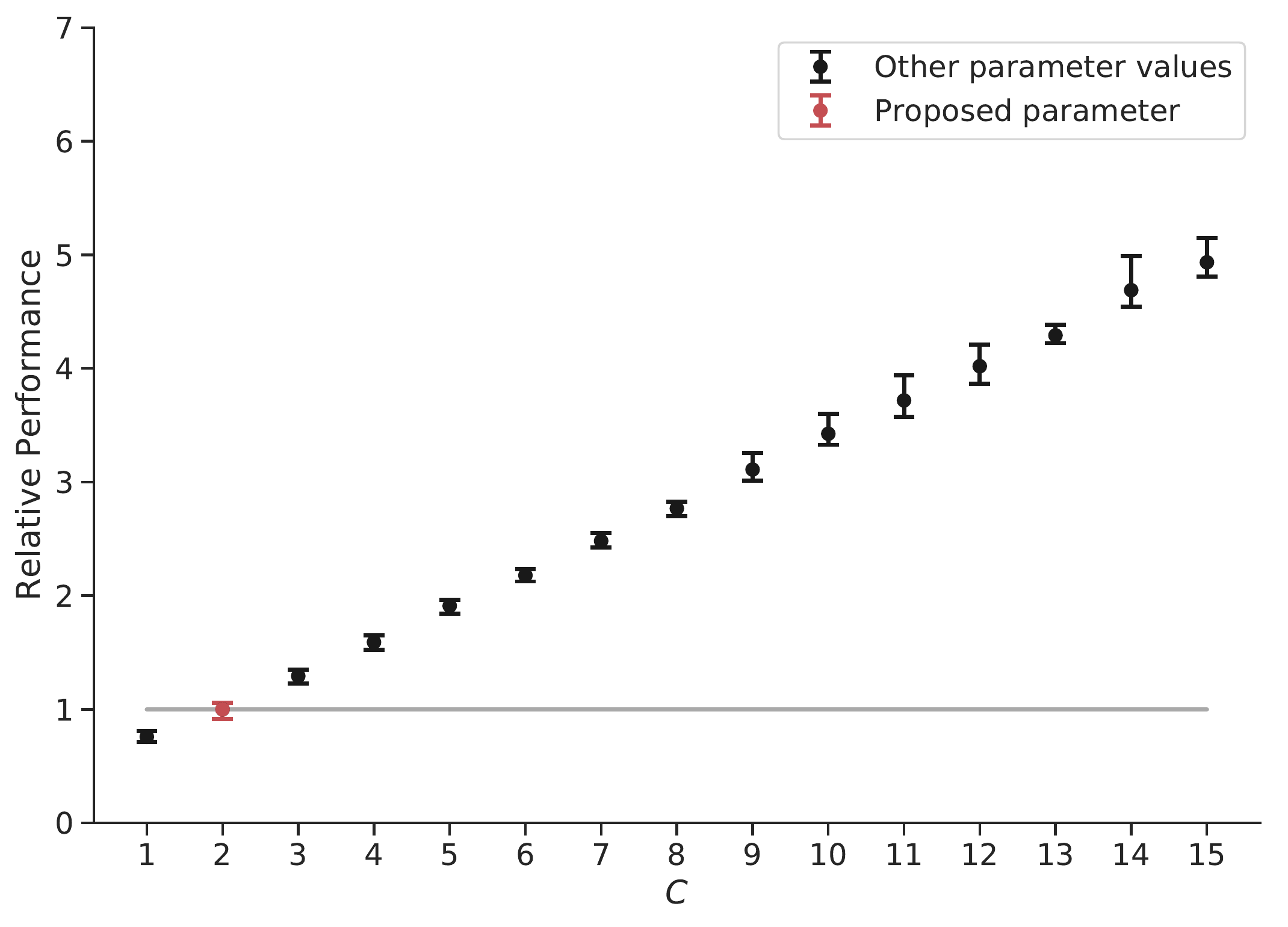}
        \caption{\texttt{\gls{lpmc}\_DC\_L}}
        \label{fig:count_DC}
    \end{subfigure}
    \hfill
    \begin{subfigure}[b]{0.49\textwidth}
        \centering
        \includegraphics[width=\textwidth]{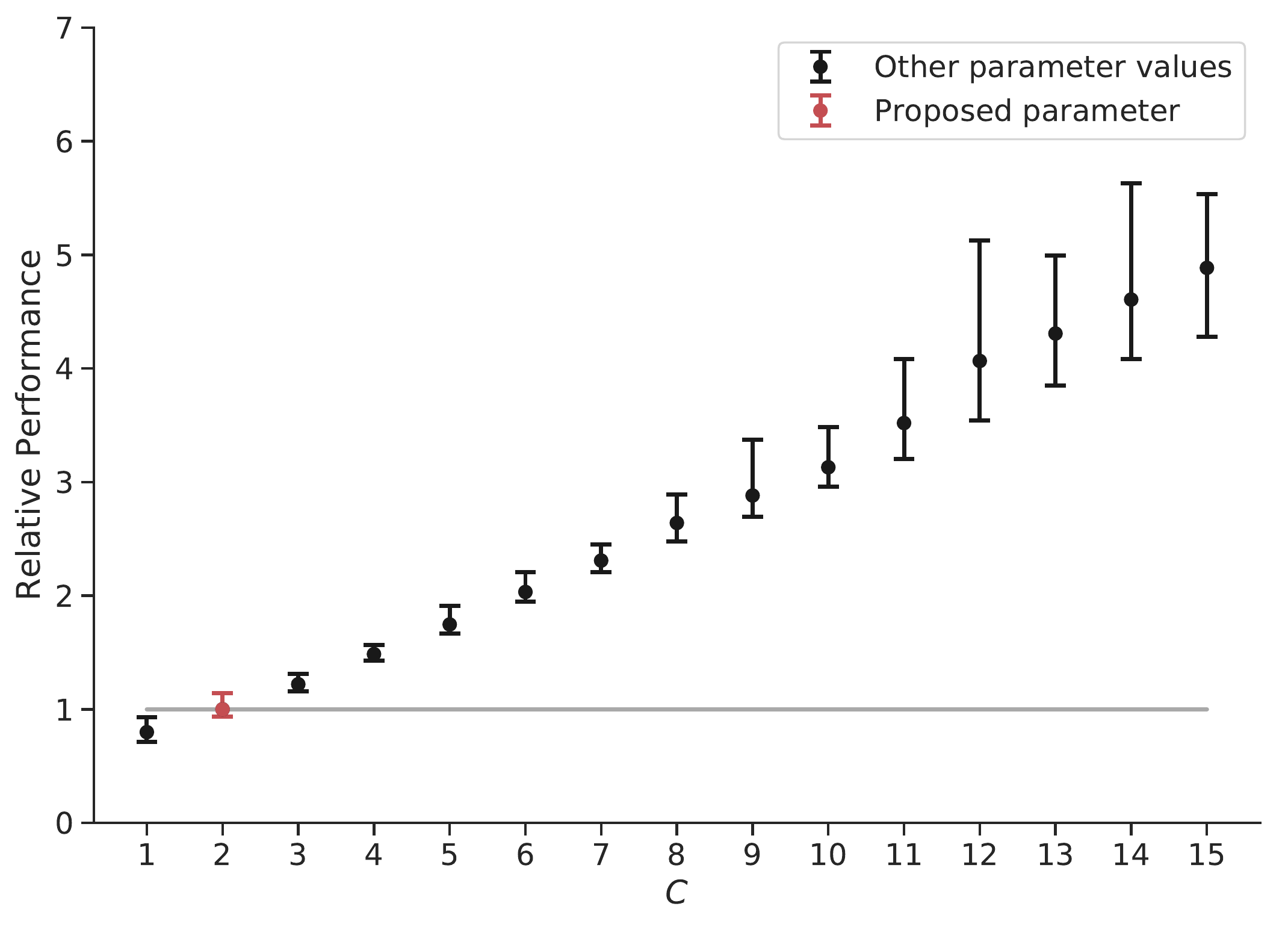}
        \caption{\texttt{\gls{lpmc}\_RR\_L}}
        \label{fig:count_RR}
    \end{subfigure}
    \hfill
    \begin{subfigure}[b]{0.49\textwidth}
        \centering
        \includegraphics[width=\textwidth]{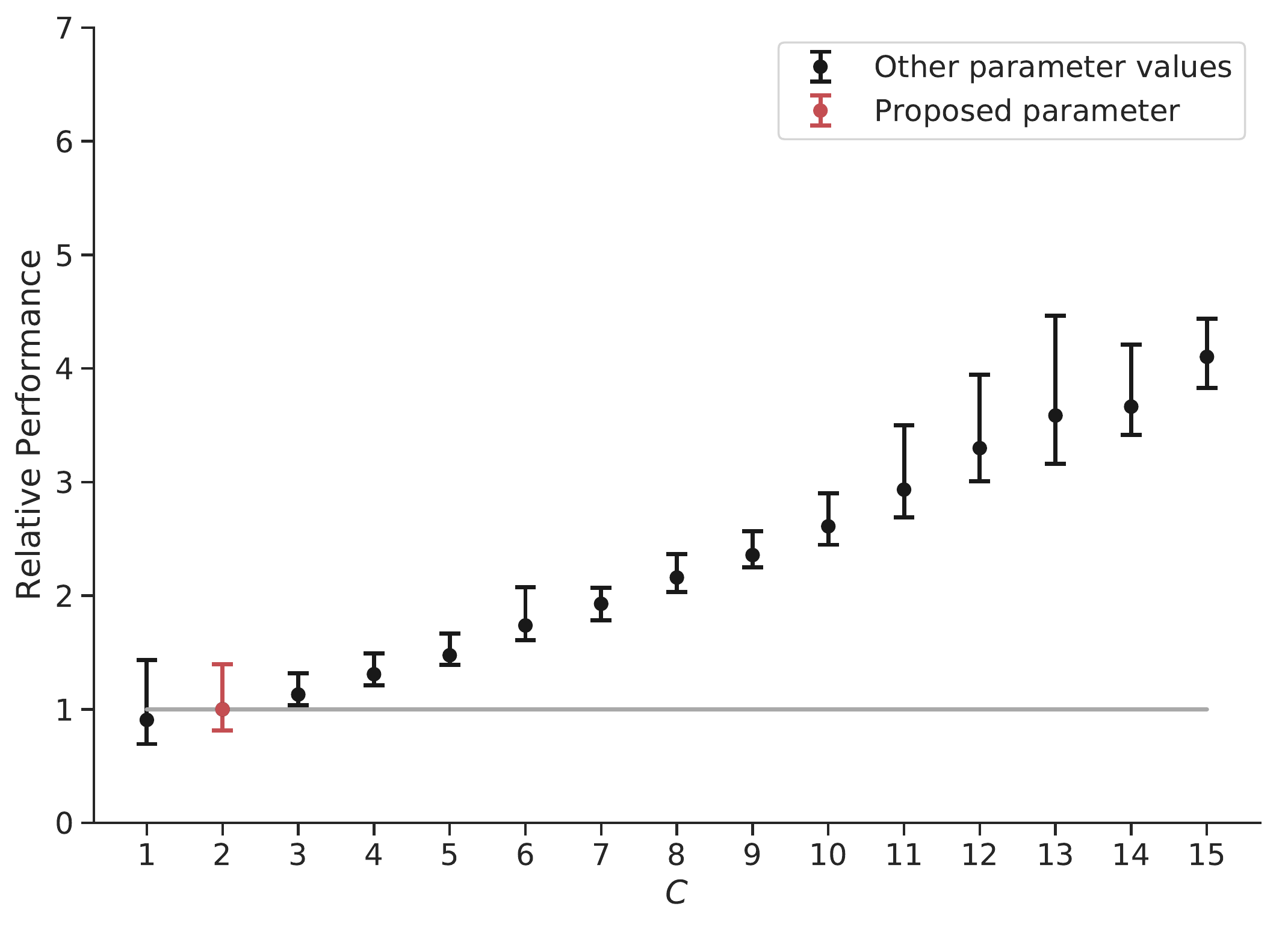}
        \caption{\texttt{\gls{lpmc}\_Full\_L}}
        \label{fig:count_Full}
    \end{subfigure}
    \hfill
    \caption{Sensitivity analysis for the parameter $C$ for the three \gls{lpmc} models using the large dataset. The black error bars correspond to the tested values and the red to the proposed value. The gray line correspond to the benchmark with the proposed parameters for the relative performance.}
    \label{fig:effect_count}
\end{figure}

Figure~\ref{fig:effect_count} shows the analysis for the parameter $C$; the maximum number of unsuccessful iterations with the same batch size. It is interesting to note that for the three models, the relationship between the execution time and this parameter value is linear. Therefore, it is evident that using a value of 1 for $C$ is faster. However, the advantage is reduced with the size of the model, as we can see when comparing Figure~\ref{fig:count_DC} and Figure~\ref{fig:count_Full}. Besides, if the count is set to 1, it could lead to triggering a false positive. Indeed, if the value of the window, $W$, is too small, an exceptionally lousy batch of data may lead to an increase of the batch size while the log likelihood can still be improved. It is thus preferable to slightly slow the execution time but increase the robustness of the algorithm. That is the reason we propose to use $C = 2$.

\begin{figure}[H]
    \centering
    \begin{subfigure}[b]{0.49\textwidth}
        \centering
        \includegraphics[width=\textwidth]{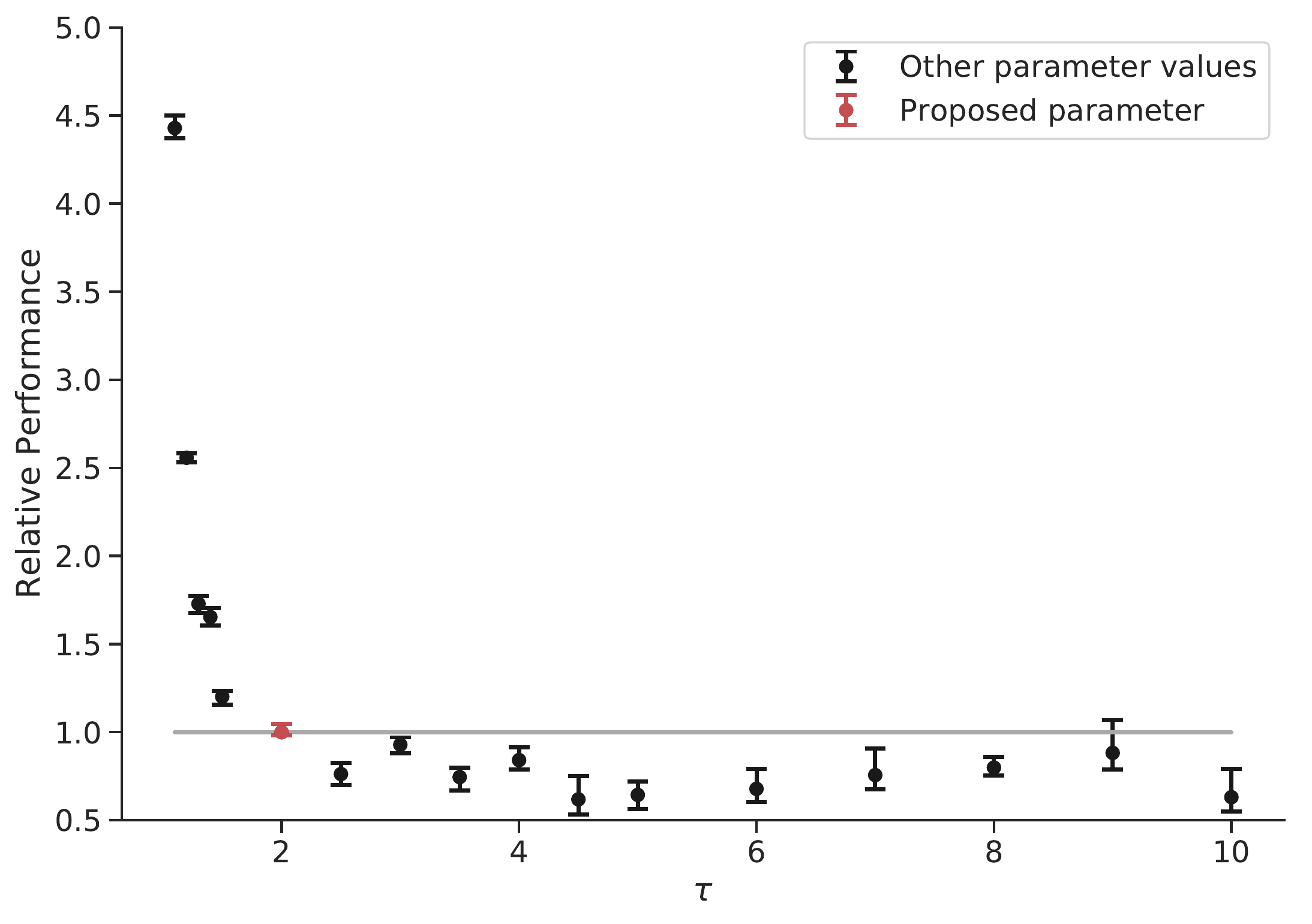}
        \caption{\texttt{\gls{lpmc}\_DC\_L}}
        \label{fig:factor_DC}
    \end{subfigure}
    \hfill
    \begin{subfigure}[b]{0.49\textwidth}
        \centering
        \includegraphics[width=\textwidth]{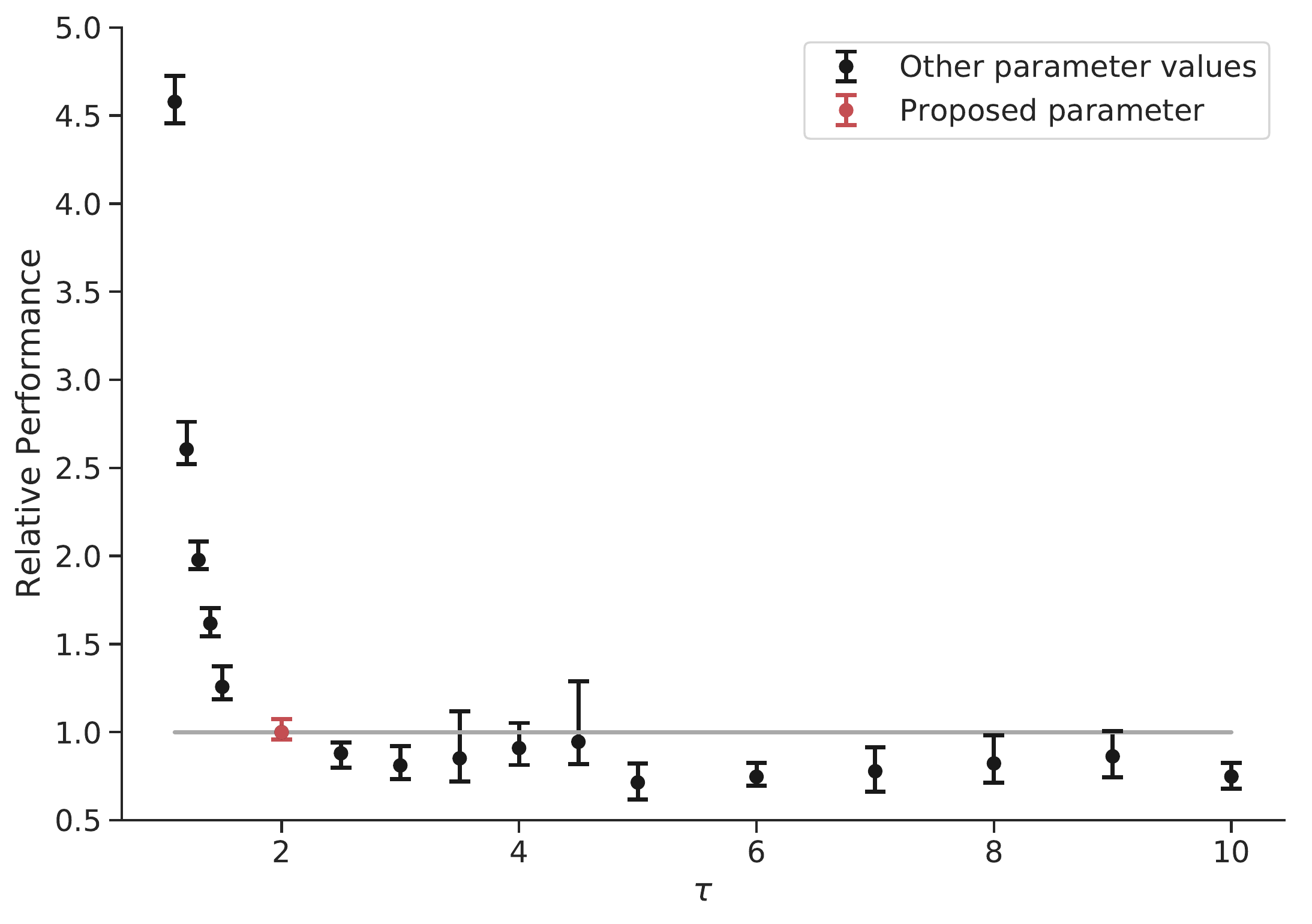}
        \caption{\texttt{\gls{lpmc}\_RR\_L}}
        \label{fig:factor_RR}
    \end{subfigure}
    \hfill
    \begin{subfigure}[b]{0.49\textwidth}
        \centering
        \includegraphics[width=\textwidth]{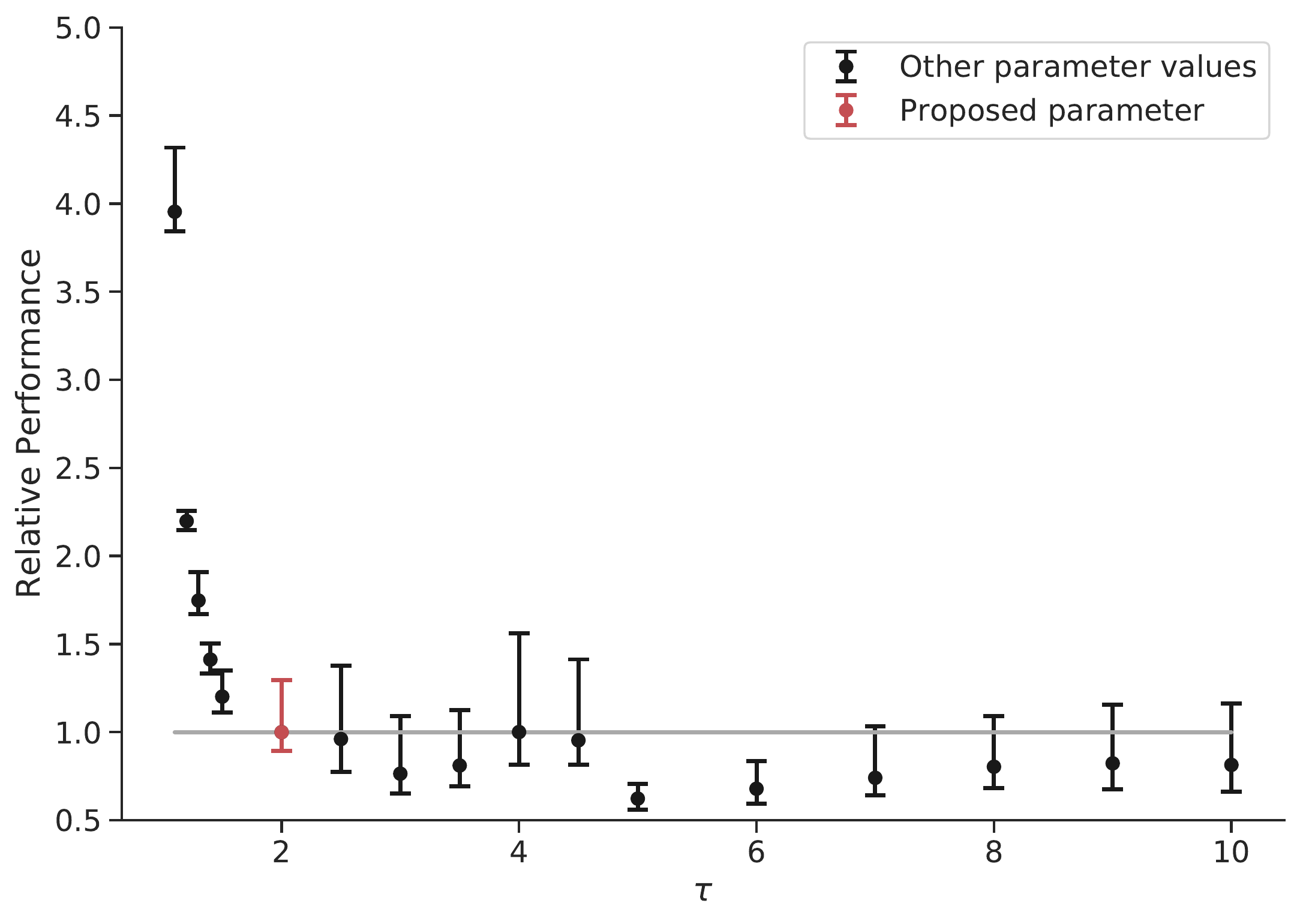}
        \caption{\texttt{\gls{lpmc}\_Full\_L}}
        \label{fig:factor_Full}
    \end{subfigure}
    \hfill
    \caption{Sensitivity analysis for the parameter $\tau$ for the three \gls{lpmc} models using the large dataset. The black error bars correspond to the tested values and the red to the proposed value. The gray line correspond to the benchmark with the proposed parameters for the relative performance.}
    \label{fig:effect_factor}
\end{figure}

Figure~\ref{fig:effect_factor} shows the analysis for the parameter $\tau$; the expansion factor for the batch size. This parameter is particularly important since it decides how fast the algorithm uses the full dataset for its iterations. As expected, a small value for $\tau$ leads to longer execution time for all three models. However, a larger value for $\tau$ is more efficient for the smaller models. Indeed, for both the \texttt{\gls{lpmc}\_DC\_L} and the \texttt{\gls{lpmc}\_RR\_L}, using larger values lead to around 20\% decrease in execution time. Since these models are quite small in the number of parameters, they already profit a lot from the starting batch size. It is then more favorable to switch the optimization algorithm and use more data. However, it seems that for the model \texttt{\gls{lpmc}\_Full\_L}, switching too soon leads to more variance in the execution time. Therefore, it is better to stay more robust and use a smaller value for $\tau$. We thus propose to use $\tau=2$.

\begin{figure}[H]
    \centering
    \begin{subfigure}[b]{0.49\textwidth}
        \centering
        \includegraphics[width=\textwidth]{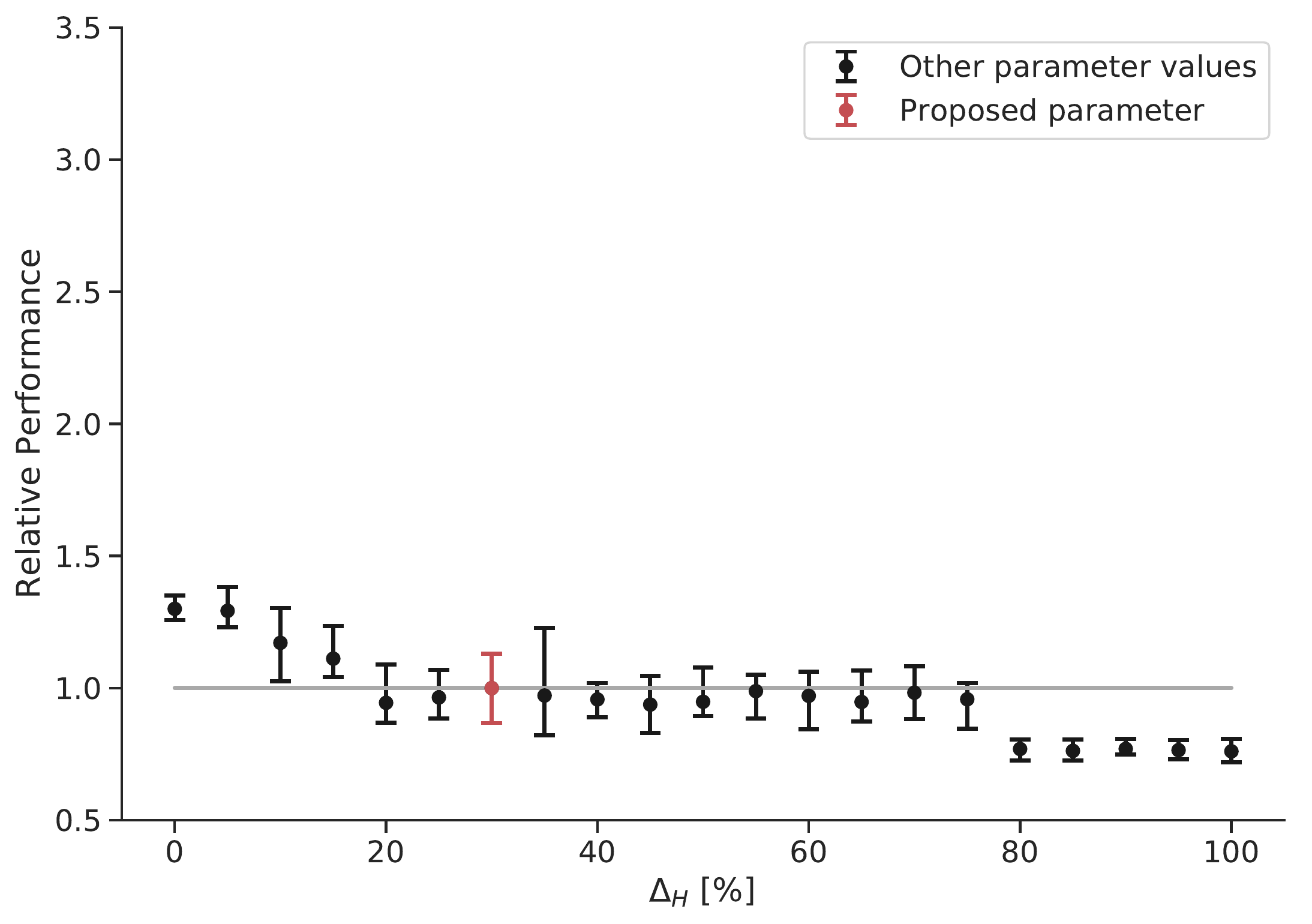}
        \caption{\texttt{\gls{lpmc}\_DC\_L}}
        \label{fig:hybrid_DC}
    \end{subfigure}
    \hfill
    \begin{subfigure}[b]{0.49\textwidth}
        \centering
        \includegraphics[width=\textwidth]{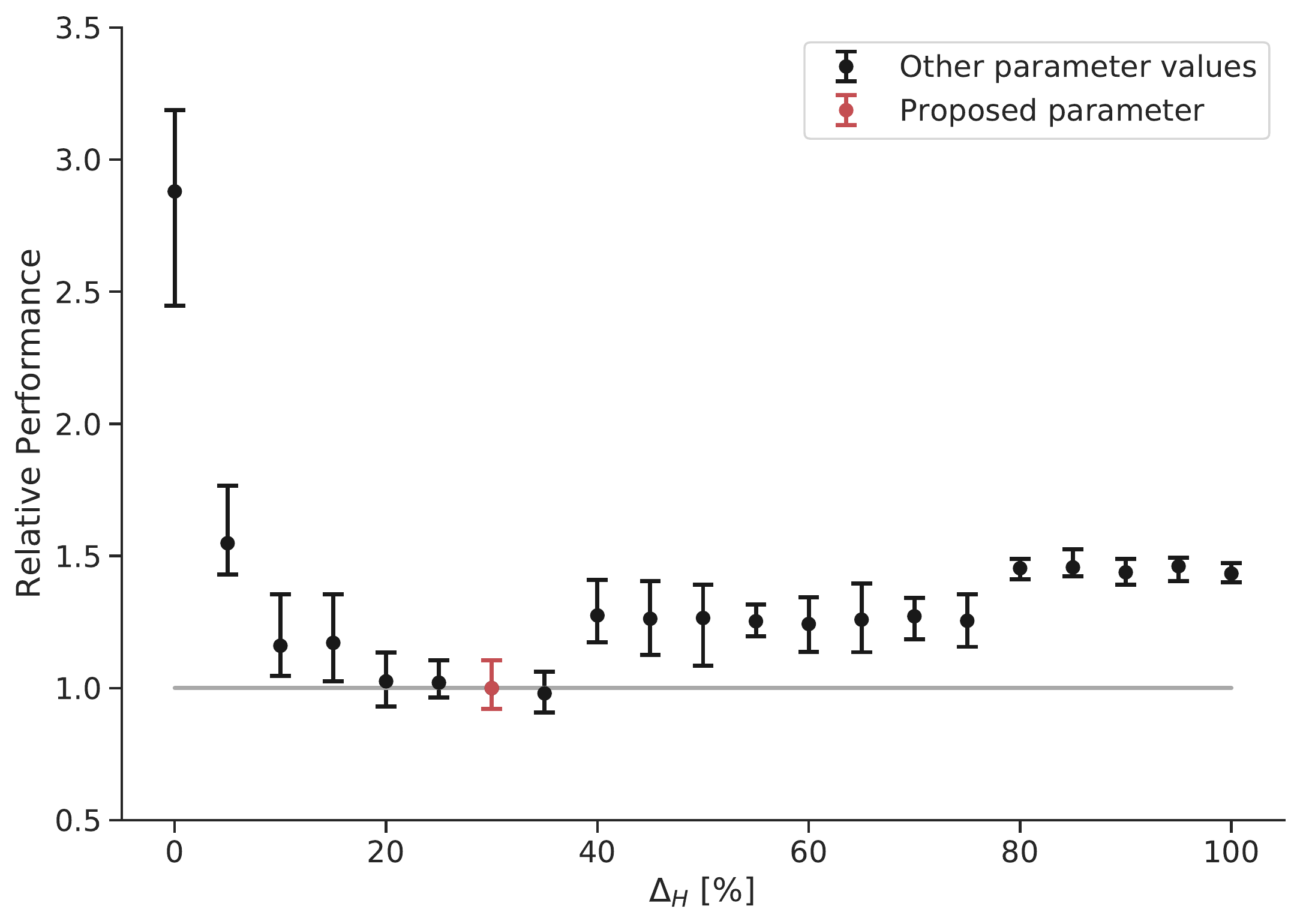}
        \caption{\texttt{\gls{lpmc}\_RR\_L}}
        \label{fig:hybrid_RR}
    \end{subfigure}
    \hfill
    \begin{subfigure}[b]{0.49\textwidth}
        \centering
        \includegraphics[width=\textwidth]{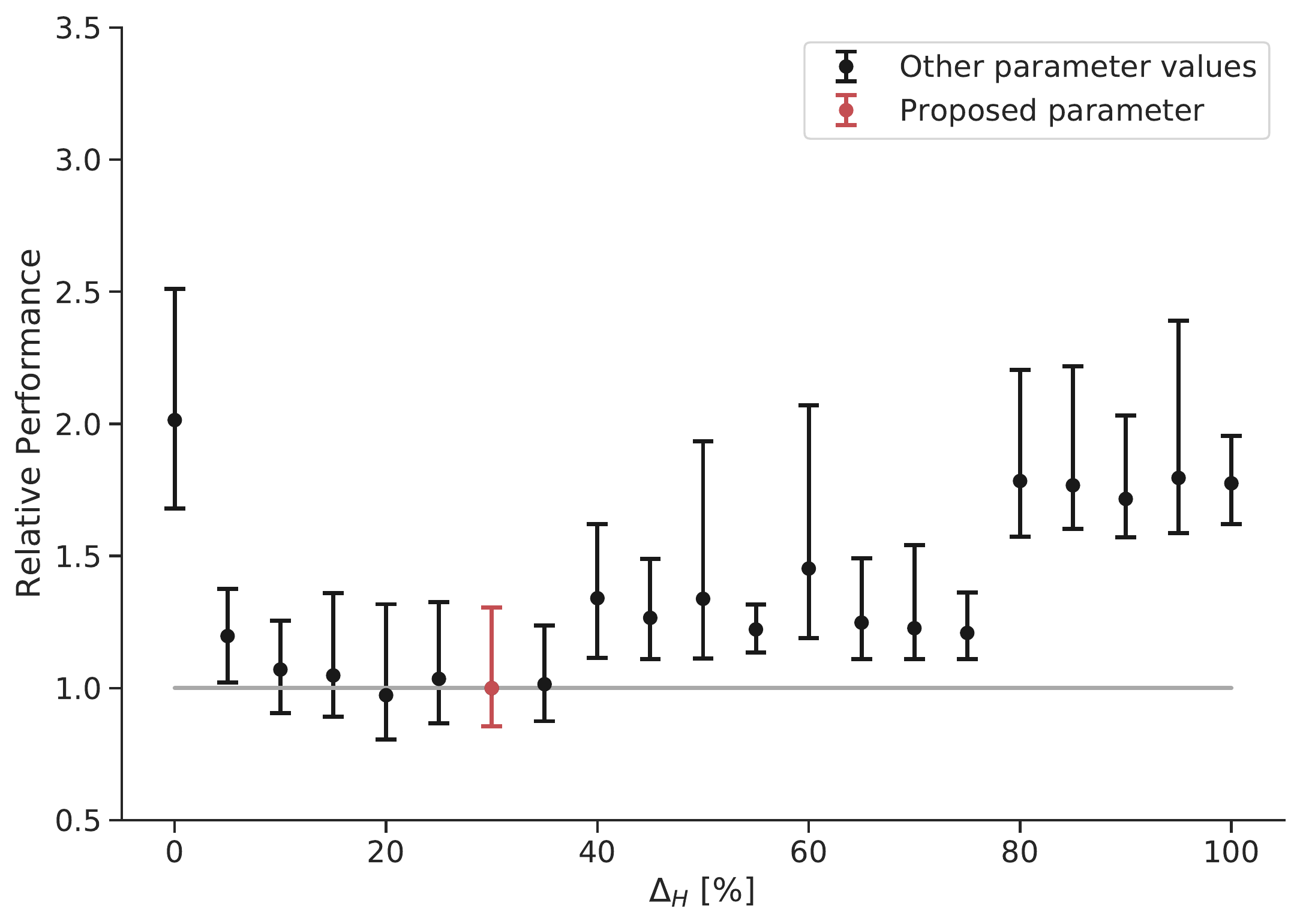}
        \caption{\texttt{\gls{lpmc}\_Full\_L}}
        \label{fig:hybrid_Full}
    \end{subfigure}
    \hfill
    \caption{Sensitivity analysis for the parameter $\Delta_H$ for the three \gls{lpmc} models using the large dataset. The black error bars correspond to the tested values and the red to the proposed value. The gray line correspond to the benchmark with the proposed parameters for the relative performance.}
    \label{fig:effect_hybrid}
\end{figure}

Figure~\ref{fig:effect_hybrid} shows the analysis for the parameter $\Delta_H$; the threshold for hybridization. Quite interestingly, this parameter does not influence much the execution time for the smallest model, the \texttt{\gls{lpmc}\_DC\_L}. As discussed in Section~\ref{sec:performance}, the \texttt{HAMABS} algorithm is not the fastest algorithm for the \texttt{\gls{lpmc}\_DC} models. Also, the execution is so small that the difference between Newton's method and BFGS is small. In addition, we see in Table~\ref{tab:results_time_1}, that Newton's method is faster than \gls{bfgs}$^{-1}$ on this particular model. Therefore, it is better to switch as late as possible. For the slightly larger model, \texttt{\gls{lpmc}\_RR\_L}, high thresholds increase the execution time. As shown in Table~\ref{tab:results_time_1}, the \gls{bfgs} methods are faster to optimize the models \texttt{\gls{lpmc}\_RR}. Therefore, it is expected that switching too late to \gls{bfgs} is increasing the execution time. However, it seems that switching as soon as possible to BFGS is the most efficient for this model. This is most likely thanks to the help of the Hessian computation at the first step. Indeed, the \gls{bfgs} algorithm generally starts with an Identity matrix. However, if we first perform a Newton step with the computation of the Hessian, it gives a good approximation as the starting point. Therefore, for the medium models, a smaller threshold for the hybridization is recommended. However, since the goal is to optimize large models as soon as possible, it is more important to use parameters specifically proposed for these models. As seen in Figure~\ref{fig:hybrid_Full}, the best switch appears at around 30\% of the data. It is slower to switch too soon since \gls{bfgs} still take more time than Newton's method to perform the early steps. And it is also slower to switch later since Newton's method takes too much time to compute the Hessian. Therefore, the proposed value is $\Delta_H=30$.

\begin{figure}[H]
    \centering
    \begin{subfigure}[b]{0.49\textwidth}
        \centering
        \includegraphics[width=\textwidth]{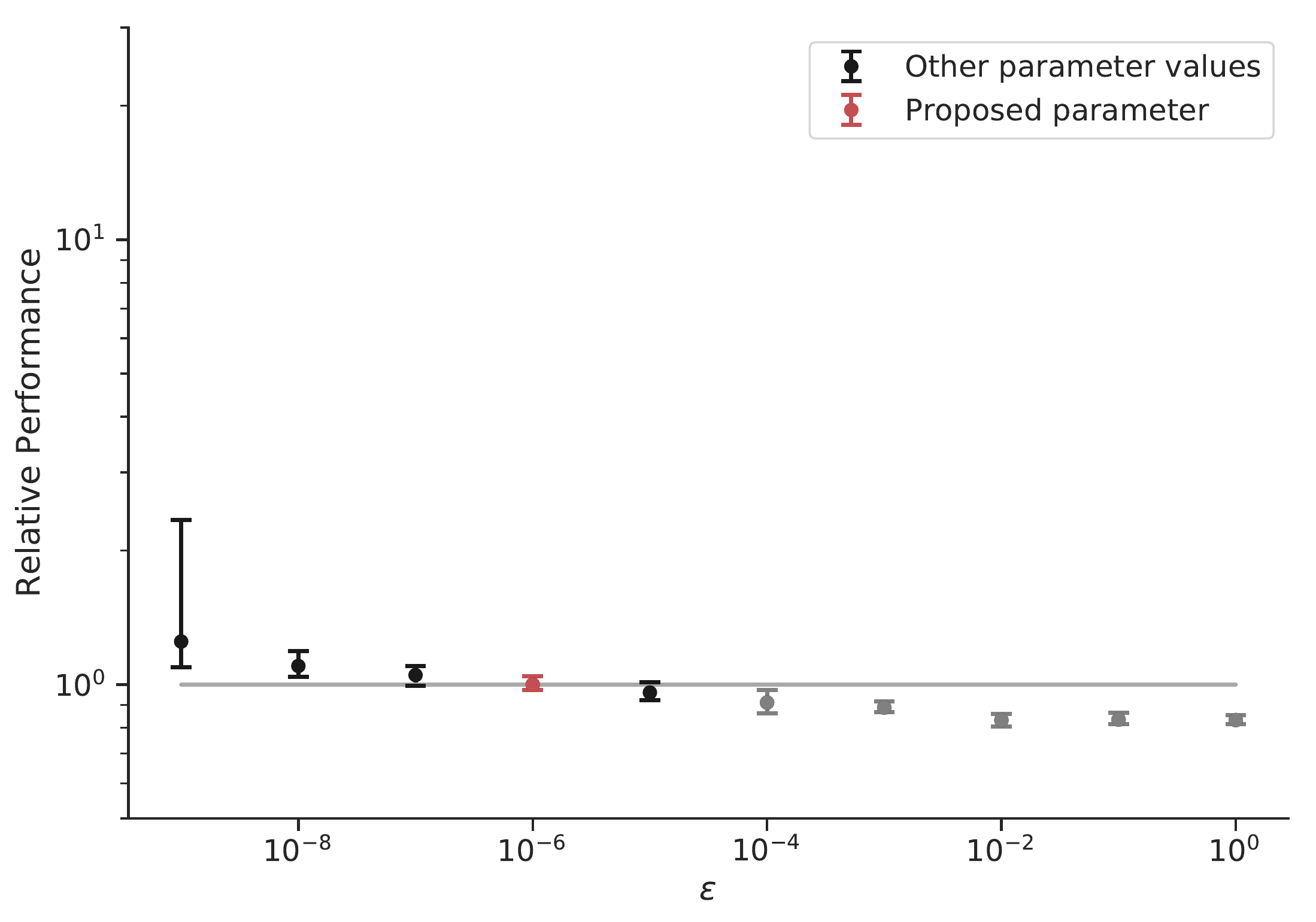}
        \caption{\texttt{\gls{lpmc}\_DC\_L}}
        \label{fig:stop_crit_DC}
    \end{subfigure}
    \hfill
    \begin{subfigure}[b]{0.49\textwidth}
        \centering
        \includegraphics[width=\textwidth]{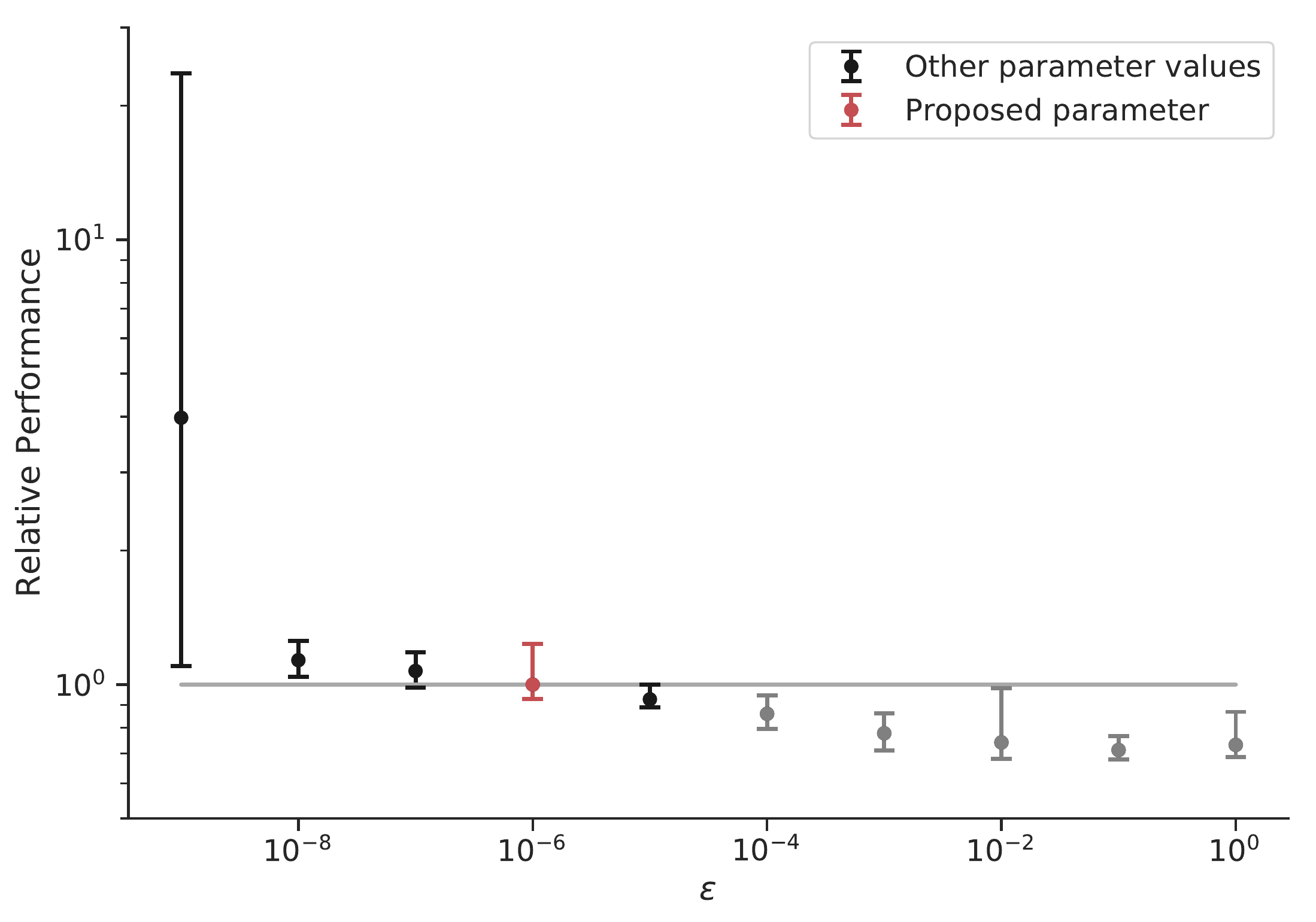}
        \caption{\texttt{\gls{lpmc}\_RR\_L}}
        \label{fig:stop_crit_RR}
    \end{subfigure}
    \hfill
    \begin{subfigure}[b]{0.49\textwidth}
        \centering
        \includegraphics[width=\textwidth]{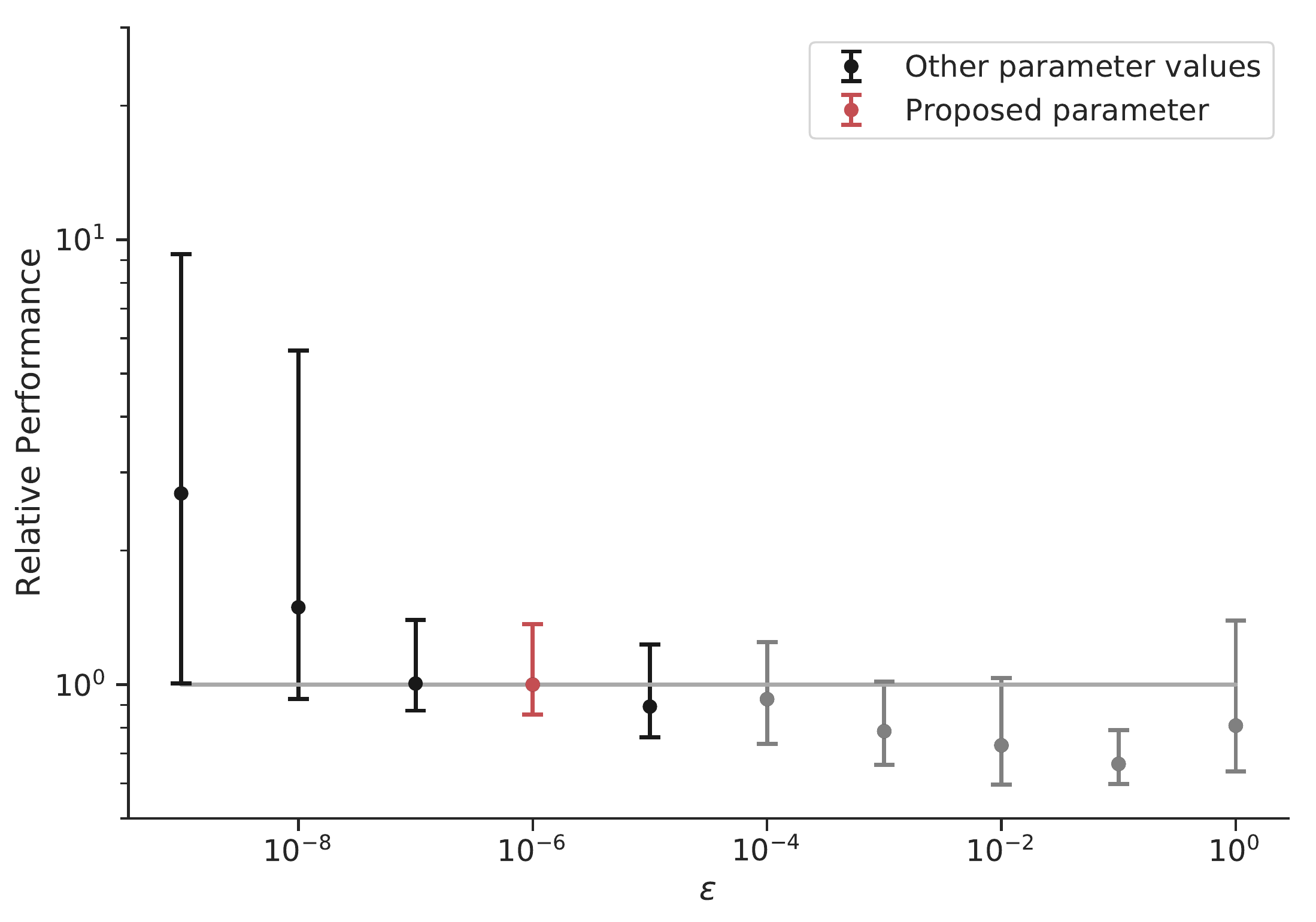}
        \caption{\texttt{\gls{lpmc}\_Full\_L}}
        \label{fig:stop_crit_Full}
    \end{subfigure}
    \hfill
    \caption{Sensitivity analysis for the parameter $\varepsilon$ for the three \gls{lpmc} models using the large dataset. The black error bars correspond to the tested values, the red to the proposed value, and the gray to optimization that were stopped too early. The gray line correspond to the benchmark with the proposed parameters for the relative performance.}
    \label{fig:effect_stop_crit}
\end{figure}

Finally, Figure~\ref{fig:effect_stop_crit} shows the analysis for the parameter $\varepsilon$; the threshold for the stopping criterion. As shown in the three figures of Figure~\ref{fig:effect_stop_crit}, using a stopping criterion that is too high is faster. However, it also leads to incorrect results. Indeed, the algorithm is stopping too soon, and the optimization has converged to the optimal point yet. We also see that using a stopping criterion too small leads to a significant increase in the execution time. In this case, it is unreasonable to try to be that precise. Therefore, a good value is a compromise between precision and speed. As shown in all these figures, a value for the stopping criterion between $10^{-8}$ and $10^{-5}$ is acceptable. We propose to use $10^{-6}$ even if other values can be used.

As seen with the sensitivity analysis, the parameters might depend on the model's size. However, the goal of this article is to be able to speed up the optimization process of large choice models. Therefore, we selected parameters that lead to an improvement on these large models. Besides, half the execution time of a small model would only result in a gain of a few seconds. On the other hand, the same speedup would lead to minutes or even hours on larger models. It is, therefore, more rewarding to speed up the larger models.

\section{Future work and conclusion}
\label{sec:conclusion}

In this article, we present three primary contributions for the estimation of \glspl{dcm}: the use of stochastic hessian, an adaptive batch size method, and the hybridization between optimization algorithms. We test 15 different algorithms ranging from standard algorithms to stochastic hybrid adaptive batch size algorithms. We prove that the use of an adaptive batch size technique is beneficial for the optimization time. Besides, since the \gls{amabs} method can be used with different optimization algorithms, we created three hybrid \gls{amabs} algorithms. We have shown that the best of these methods is the \texttt{HAMABS} algorithm. It speeds up the optimization by a factor of 23 on the \texttt{\gls{mtmc}} model containing 247 parameters. Therefore, the use of faster algorithms opens the research to new possibilities for the future of choice modeling. As a concrete example, faster optimization time allows researchers to test many more specifications in the same amount of time. This can thus be used to develop automatic utility specification techniques to speed up the modelization of \glspl{dcm}.
\\\\
For the future, we would like to work on two different improvements. The first one concerns the hybridization. The current way of doing it depends on the starting batch size. Indeed, due to the geometrical rule to increase the batch size in \texttt{AMABS}, it would be possible to miss the 30\% mark to switch the optimization algorithm. Therefore, we would like to work on a better switch for the hybridization. One possible direction is to compare the improvement made by each algorithm in one step over the time it takes to do it. We would then have a metric in the percentage of improvement over seconds, and we could easily decide to switch the algorithm when BFGS leads to more improvements per second. The second improvement can be made on the rule for the update of the batch size in \texttt{AMABS}. Indeed, the geometrical rule seems to work well. However, it is possible that combining multiple rules could lead to faster optimization time. The final improvement is to integrate algorithms dealing with bounds. Indeed, nested logit and MEV models require bounds for some parameters. It is thus essential to integrate them into future optimization algorithms. Finally, we would like to implement the final version of this algorithm in Pandas Biogeme.

\newpage
\section{Appendix}
\label{sec:appendix}

\subsection{Execution time - Tables}
\label{sec:appendix_perf_time}

Tables~\ref{tab:results_time_1} and \ref{tab:results_time_2} show the average estimation time with the standard deviation for each models optimized by each algorithm. The bold values with the gray background represent the fastest algorithms for each model. The gray values mean that the algorithms were not able to converge in the required number of epochs (1000).

\subsection{Number of epochs - Tables}
\label{sec:appendix_perf_epochs}

Tables~\ref{tab:results_epochs_1} and \ref{tab:results_epochs_2} show the average number of epochs with the standard deviation used by each algorithm to optimize each model. The bold values with the gray background represent the algorithms that have used the least epochs for each model. The gray values mean that the algorithms were not able to converge in the required number of epochs (1000).

\begin{landscape}
\vspace*{\fill}
\begin{table}[H]
  \renewcommand\arraystretch{1.1}
  \centering
  \begin{tabular}{l||c|c|c|c|c|c}
    \multirow{2}{*}{Algorithms} & \multicolumn{6}{c}{Models} \\ \cline{2-7}
  & \texttt{LPMC\_DC\_S} & \texttt{LPMC\_DC\_M} & \texttt{LPMC\_DC\_L} & \texttt{LPMC\_RR\_S} & \texttt{LPMC\_RR\_M} & \texttt{LPMC\_RR\_L} \\ \hline\hline
    \texttt{GD} & \color{gray}$52.16 \pm 0.24$ & \color{gray}$95.07 \pm 0.45$ & \color{gray}$140.95 \pm 0.49$ & \color{gray}$303.21 \pm 0.25$ & \color{gray}$567.61 \pm 0.45$ & \color{gray}$820.38 \pm 1.46$ \\\hline
    \texttt{BFGS} & $6.78 \pm 0.03$ & $13.11 \pm 0.06$ & $17.96 \pm 0.06$ & $183.83 \pm 0.50$ & $337.25 \pm 0.30$ & $492.90 \pm 0.39$ \\\hline
    \texttt{BFGS$^{-1}$} & $7.20 \pm 0.02$ & $12.96 \pm 0.10$ & $19.61 \pm 0.12$ & $177.66 \pm 0.25$ & $332.93 \pm 0.40$ & $473.32 \pm 0.98$ \\\hline
    \texttt{TR-BFGS} & $5.40 \pm 0.04$ & $11.05 \pm 0.11$ & $17.17 \pm 0.07$ & $227.83 \pm 0.51$ & $405.82 \pm 0.60$ & \color{gray}$627.43 \pm 0.67$ \\\hline
    \texttt{NM} & $0.65 \pm 0.01$ & $1.23 \pm 0.01$ & $1.95 \pm 0.02$ & $36.05 \pm 0.98$ & $70.28 \pm 1.50$ & $126.85 \pm 2.18$ \\\hline
    \texttt{TR} & \cellcolor{black!15}\boldmath$0.40 \pm 0.01$ & \cellcolor{black!15}\boldmath$0.74 \pm 0.01$ & \cellcolor{black!15}\boldmath$1.03 \pm 0.01$ & $23.63 \pm 0.49$ & $46.16 \pm 0.66$ & $67.13 \pm 1.02$ \\\hline\hline
    \texttt{GD-ABS} & \color{gray}$50.53 \pm 0.32$ & \color{gray}$92.80 \pm 0.51$ & \color{gray}$138.20 \pm 0.52$ & \color{gray}$310.59 \pm 1.21$ & \color{gray}$568.12 \pm 1.74$ & \color{gray}$824.15 \pm 2.54$ \\\hline
    \texttt{BFGS-ABS} & $6.46 \pm 0.29$ & $10.80 \pm 0.37$ & $15.27 \pm 0.63$ & $169.83 \pm 1.58$ & $314.56 \pm 3.83$ & $474.42 \pm 5.53$ \\\hline
    \texttt{BFGS$^{-1}$-ABS} & $6.85 \pm 0.10$ & $11.13 \pm 0.14$ & $15.62 \pm 0.16$ & $169.36 \pm 0.62$ & $312.78 \pm 2.79$ & $460.49 \pm 1.92$ \\\hline
    \texttt{TR-BFGS-ABS} & $6.90 \pm 0.37$ & $12.49 \pm 0.76$ & $21.56 \pm 1.04$ & $194.49 \pm 8.57$ & $387.56 \pm 12.13$ & $590.62 \pm 18.80$ \\\hline
    \texttt{NM-ABS} & $7.64 \pm 18.97$ & $36.79 \pm 52.75$ & $116.95 \pm 74.20$ & $29.72 \pm 2.41$ & $53.28 \pm 3.10$ & $76.44 \pm 5.78$ \\\hline
    \texttt{TR-ABS} & $3.20 \pm 0.06$ & $6.43 \pm 0.15$ & $11.50 \pm 0.13$ & $50.64 \pm 2.10$ & $97.66 \pm 4.82$ & $175.75 \pm 5.21$ \\\hline\hline
    \texttt{H-NM-ABS} & $2.83 \pm 0.15$ & $4.77 \pm 0.25$ & $6.97 \pm 0.39$ & $29.80 \pm 1.38$ & $48.46 \pm 2.50$ & $20.78 \pm 1.96$ \\\hline
    \texttt{H-TR-ABS} & $2.16 \pm 0.14$ & $3.82 \pm 0.20$ & $6.88 \pm 0.38$ & $34.97 \pm 6.11$ & $60.87 \pm 13.06$ & $157.69 \pm 10.85$ \\\hline
    \texttt{HAMABS} & $1.86 \pm 0.12$ & $3.11 \pm 0.20$ & $4.59 \pm 0.32$ & \cellcolor{black!15}\boldmath$11.98 \pm 1.23$ & \cellcolor{black!15}\boldmath$18.46 \pm 1.06$ & \cellcolor{black!15}\boldmath$18.14 \pm 1.06$ \\
  \end{tabular}
  \caption{Time in seconds used for the estimation of the models \texttt{LPMC\_DC} and \texttt{LPMC\_RR} by all the algorithms presented in Table~\ref{tab:algorithms}. The values in light gray mean that the algorithms was not able to converge in the required number of eqpochs. The values in bold, in a gray cell, correspond to the the fastest optimization time.}
  \label{tab:results_time_1}
\end{table}
\vspace*{\fill}
\end{landscape}

\begin{landscape}
\vspace*{\fill}
\begin{table}[H]
  \renewcommand\arraystretch{1.1}
  \centering
  \begin{tabular}{l||c|c|c||c}
    \multirow{2}{*}{Algorithms} & \multicolumn{4}{c}{Models} \\ \cline{2-5}
  & \texttt{LPMC\_Full\_S} & \texttt{LPMC\_Full\_M} & \texttt{LPMC\_Full\_L} & \texttt{MTMC} \\ \hline\hline
    \texttt{GD} & \color{gray}$2785.66 \pm 28.12$ & \color{gray}$5344.42 \pm 33.60$ & \color{gray}$8031.02 \pm 49.62$ & \color{gray}$7842.92 \pm 53.75$ \\\hline
    \texttt{BFGS} & $3173.86 \pm 22.37$ & $6289.47 \pm 34.27$ & $9333.66 \pm 63.08$ & \color{gray}$10308.49 \pm 61.36$ \\\hline
    \texttt{BFGS$^{-1}$} & $3129.08 \pm 10.83$ & $5929.40 \pm 27.85$ & $8812.09 \pm 66.60$ & \color{gray}$10090.74 \pm 56.43$ \\\hline
    \texttt{TR-BFGS} & \color{gray}$2043.68 \pm 16.10$ & \color{gray}$4014.60 \pm 24.31$ & \color{gray}$5861.06 \pm 39.09$ & \color{gray}$5568.99 \pm 38.86$ \\\hline
    \texttt{NM} & $1398.81 \pm 52.16$ & $3230.83 \pm 69.75$ & $3984.04 \pm 91.00$ & $17199.76 \pm 190.18$ \\\hline
    \texttt{TR} & $540.07 \pm 24.73$ & $1021.88 \pm 38.01$ & $1501.90 \pm 53.11$ & $10613.62 \pm 152.59$ \\\hline\hline
    \texttt{GD-ABS} & \color{gray}$2782.30 \pm 23.46$ & \color{gray}$5508.52 \pm 27.21$ & \color{gray}$7989.58 \pm 64.22$ & \color{gray}$7984.13 \pm 65.41$ \\\hline
    \texttt{BFGS-ABS} & $3163.19 \pm 57.56$ & $6119.99 \pm 109.08$ & $8944.95 \pm 130.30$ & \color{gray}$10225.68 \pm 59.41$ \\\hline
    \texttt{BFGS$^{-1}$-ABS} & $3021.65 \pm 28.72$ & $5814.22 \pm 47.75$ & $8721.32 \pm 92.05$ & \color{gray}$10083.01 \pm 58.33$ \\\hline
    \texttt{TR-BFGS-ABS} & \color{gray}$2035.42 \pm 17.94$ & \color{gray}$4001.22 \pm 29.88$ & \color{gray}$5796.55 \pm 47.00$ & \color{gray}$5586.36 \pm 39.60$ \\\hline
    \texttt{NM-ABS} & $4359.60 \pm 14019.82$ & $1982.06 \pm 229.38$ & $2721.13 \pm 270.25$ & $14535.94 \pm 501.67$ \\\hline
    \texttt{TR-ABS} & $1203.69 \pm 84.31$ & $2240.94 \pm 99.74$ & $4089.42 \pm 134.02$ & $13695.15 \pm 757.32$ \\\hline\hline
    \texttt{H-NM-ABS} & $522.65 \pm 47.01$ & $939.05 \pm 85.38$ & $1376.32 \pm 93.92$ & $2749.46 \pm 103.71$ \\\hline
    \texttt{H-TR-ABS} & $591.33 \pm 80.36$ & $1085.78 \pm 137.27$ & $2209.45 \pm 160.07$ & $7633.65 \pm 386.77$ \\\hline
    \texttt{HAMABS} & \cellcolor{black!15}\boldmath$257.02 \pm 42.85$ & \cellcolor{black!15}\boldmath$405.43 \pm 43.77$ & \cellcolor{black!15}\boldmath$486.31 \pm 63.38$ & \cellcolor{black!15}\boldmath$1243.95 \pm 56.21$ \\
  \end{tabular}
  \caption{Time in seconds used for the estimation of the models \texttt{LPMC\_Full} and \texttt{MTMC} by all the algorithms presented in Table~\ref{tab:algorithms}. The values in light gray mean that the algorithms was not able to converge in the required number of eqpochs. The values in bold, in a gray cell, correspond to the the fastest optimization time.}
  \label{tab:results_time_2}
\end{table}
\vspace*{\fill}
\end{landscape}

\begin{landscape}
\vspace*{\fill}
\begin{table}[H]
  \renewcommand\arraystretch{1.1}
  \centering
  \begin{tabular}{l||c|c|c|c|c|c}
    \multirow{2}{*}{Algorithms} & \multicolumn{6}{c}{Models} \\ \cline{2-7}
  & \texttt{LPMC\_DC\_S} & \texttt{LPMC\_DC\_M} & \texttt{LPMC\_DC\_L} & \texttt{LPMC\_RR\_S} & \texttt{LPMC\_RR\_M} & \texttt{LPMC\_RR\_L} \\ \hline\hline
    \texttt{GD} & \color{gray}$1000$ & \color{gray}$1000$ & \color{gray}$1000$ & \color{gray}$1000$ & \color{gray}$1000$ & \color{gray}$1000$ \\\hline
    \texttt{BFGS} & $108$ & $109$ & $99$ & $480$ & $461$ & $478$ \\\hline
    \texttt{BFGS$^{-1}$} & $111$ & $112$ & $111$ & $468$ & $464$ & $462$ \\\hline
    \texttt{TR-BFGS} & $132$ & $148$ & $151$ & $989$ & $935$ & \color{gray}$1000$ \\\hline
    \texttt{NM} & $9$ & $10$ & $11$ & $12$ & $12$ & $15$ \\\hline
    \texttt{TR} & \cellcolor{black!15}\boldmath$8$ & \cellcolor{black!15}\boldmath$8$ & \cellcolor{black!15}\boldmath$8$ & \cellcolor{black!15}\boldmath$8$ & \cellcolor{black!15}\boldmath$8$ & \cellcolor{black!15}\boldmath$8$ \\\hline\hline
    \texttt{GD-ABS} & \color{gray}$1000.71 \pm 0.09$ & \color{gray}$1000.53 \pm 0.04$ & \color{gray}$1000.28 \pm 0.02$ & \color{gray}$1000.62 \pm 0.26$ & \color{gray}$1000.57 \pm 0.16$ & \color{gray}$1000.29 \pm 0.11$ \\\hline
    \texttt{BFGS-ABS} & $85.05 \pm 5.75$ & $81.91 \pm 3.91$ & $78.88 \pm 4.63$ & $445.97 \pm 4.61$ & $442.97 \pm 6.35$ & $441.12 \pm 5.45$ \\\hline
    \texttt{BFGS$^{-1}$-ABS} & $90.02 \pm 0.97$ & $87.79 \pm 1.55$ & $85.64 \pm 1.21$ & $445.21 \pm 1.64$ & $439.54 \pm 3.81$ & $437.46 \pm 1.80$ \\\hline
    \texttt{TR-BFGS-ABS} & $101.91 \pm 12.76$ & $104.98 \pm 11.82$ & $109.02 \pm 11.43$ & $838.64 \pm 38.49$ & $880.85 \pm 28.64$ & $937.42 \pm 31.23$ \\\hline
    \texttt{NM-ABS} & $106.97 \pm 297.92$ & $304.96 \pm 455.35$ & $702.37 \pm 455.05$ & $10.13 \pm 0.78$ & $9.28 \pm 0.60$ & $9.11 \pm 0.75$ \\\hline
    \texttt{TR-ABS} & $15.62 \pm 0.34$ & $15.91 \pm 0.33$ & $20.41 \pm 0.18$ & $16.77 \pm 0.35$ & $16.71 \pm 0.59$ & $20.54 \pm 0.42$ \\\hline\hline
    \texttt{H-NM-ABS} & $30.68 \pm 2.86$ & $30.26 \pm 2.39$ & $28.07 \pm 2.40$ & $67.84 \pm 4.32$ & $60.45 \pm 4.39$ & $11.96 \pm 1.15$ \\\hline
    \texttt{H-TR-ABS} & $37.53 \pm 3.55$ & $37.54 \pm 2.78$ & $47.26 \pm 3.90$ & $147.98 \pm 30.61$ & $144.95 \pm 36.72$ & $218.49 \pm 21.42$ \\\hline
    \texttt{HAMABS} & $13.79 \pm 1.70$ & $13.50 \pm 2.05$ & $12.57 \pm 1.93$ & $15.31 \pm 2.53$ & $13.95 \pm 1.67$ & $9.55 \pm 0.56$ \\
  \end{tabular}
  \caption{Number of epochs used for the estimation of the models \texttt{LPMC\_DC} and \texttt{LPMC\_RR} by all the algorithms presented in Table~\ref{tab:algorithms}. The values in light gray mean that the algorithms was not able to converge in the required number of eqpochs. The values in bold, in a gray cell, correspond to the the fastest optimization time.}
  \label{tab:results_epochs_1}
\end{table}
\vspace*{\fill}
\end{landscape}

\begin{landscape}
\vspace*{\fill}
\begin{table}[H]
  \renewcommand\arraystretch{1.1}
  \centering
  \begin{tabular}{l||c|c|c||c}
    \multirow{2}{*}{Algorithms} & \multicolumn{4}{c}{Models} \\ \cline{2-5}
  & \texttt{LPMC\_Full\_S} & \texttt{LPMC\_Full\_M} & \texttt{LPMC\_Full\_L} & \texttt{MTMC} \\ \hline\hline
    \texttt{GD} & \color{gray}$1000$ & \color{gray}$1000$ & \color{gray}$1000$ & \color{gray}$1000$ \\\hline
    \texttt{BFGS} & $872$ & $901$ & $885$ & \color{gray}$1000$ \\\hline
    \texttt{BFGS$^{-1}$} & $878$ & $859$ & $868$ & \color{gray}$1000$ \\\hline
    \texttt{TR-BFGS} & \color{gray}$1000$ & \color{gray}$1000$ & \color{gray}$1000$ & \color{gray}$1000$ \\\hline
    \texttt{NM} & $20$ & $24$ & $20$ & $23$ \\\hline
    \texttt{TR} & \cellcolor{black!15}\boldmath$7$ & \cellcolor{black!15}\boldmath$7$ & \cellcolor{black!15}\boldmath$7$ & \cellcolor{black!15}\boldmath$14$ \\\hline\hline
    \texttt{GD-ABS} & \color{gray}$1000.71 \pm 0.15$ & \color{gray}$1000.54 \pm 0.05$ & \color{gray}$1000.29 \pm 0.03$ & \color{gray}$1000.30 \pm 0.03$ \\\hline
    \texttt{BFGS-ABS} & $877.97 \pm 16.90$ & $870.00 \pm 13.58$ & $867.96 \pm 11.58$ & \color{gray}$1000.32 \pm 0.05$ \\\hline
    \texttt{BFGS$^{-1}$-ABS} & $856.91 \pm 4.67$ & $841.77 \pm 6.28$ & $843.83 \pm 3.78$ & \color{gray}$1000.37 \pm 0.06$ \\\hline
    \texttt{TR-BFGS-ABS} & \color{gray}$1000.67 \pm 0.07$ & \color{gray}$1000.28 \pm 0.04$ & \color{gray}$1000.79 \pm 0.37$ & \color{gray}$1000.74 \pm 0.30$ \\\hline
    \texttt{NM-ABS} & $65.91 \pm 214.46$ & $15.21 \pm 1.62$ & $14.11 \pm 1.20$ & $20.12 \pm 0.60$ \\\hline
    \texttt{TR-ABS} & $17.41 \pm 1.06$ & $17.10 \pm 0.50$ & $21.34 \pm 0.39$ & $18.40 \pm 1.03$ \\\hline\hline
    \texttt{H-NM-ABS} & $98.20 \pm 4.60$ & $94.75 \pm 5.36$ & $99.43 \pm 4.52$ & $145.89 \pm 8.99$ \\\hline
    \texttt{H-TR-ABS} & $219.83 \pm 37.59$ & $216.80 \pm 34.95$ & $344.72 \pm 27.79$ & $564.62 \pm 81.41$ \\\hline
    \texttt{HAMABS} & $24.98 \pm 2.16$ & $20.42 \pm 1.84$ & $21.36 \pm 2.05$ & $18.63 \pm 1.58$ \\
  \end{tabular}
  \caption{Number of epochs used for the estimation of the models \texttt{LPMC\_Full} and \texttt{MTMC} by all the algorithms presented in Table~\ref{tab:algorithms}. The values in light gray mean that the algorithms was not able to converge in the required number of eqpochs. The values in bold, in a gray cell, correspond to the the fastest optimization time.}
  \label{tab:results_epochs_2}
\end{table}
\vspace*{\fill}
\end{landscape}

\newpage
\bibliography{LedLurHilBie_OptiDCM}

\end{document}